\documentclass[leqno]{elsart1p}

\usepackage{amssymb}

\usepackage{amsmath}

\usepackage{color}

\newtheorem{theorem}{Theorem}[section]
\newtheorem{lemma}[theorem]{Lemma}

\newtheorem{corollary}[theorem]{Corollary}

\newtheorem{proposition}[theorem]{Proposition}
\newtheorem{remark}[theorem]{Remark}

\newcommand{\R}{\mathbb{R}}
\newcommand{\C}{{\mathbb{C}}}
\newcommand{\N}{{\mathbb{N}}}
\renewcommand{\j}{{\mathcal{J}}}
\newcommand{\eps}{\varepsilon}
\newcommand{\m}{{\mathcal{M}}}
\newcommand{\M}{{\mathcal{M}}}
\renewcommand{\S}{{\mathcal{S}}}
\renewcommand{\L}{{\mathcal{L}}}
\newcommand{\equ}[1]{(\ref{#1})}
\newcommand{\opm}[1]{\M[{#1}]}
\newcommand{\ops}[1]{\S[{#1}]}
\newcommand{\opl}[1]{\L[{#1}]}
\newcommand{\dem}[1]{\medskip \noindent {\bf #1}}
\newcommand{\fdem}{\hfill$\square$\medskip}
\newcommand{\under}[1]{\underline{#1}}
\newcommand{\nli}[2]{\|{#1}\|_{L^{\infty}{(#2)}}}
\renewcommand{\O}{{\Omega}}
\renewcommand{\o}{{\omega}}

\numberwithin{equation}{section}

\begin{document}

\begin{frontmatter}



\title{Non-local anisotropic dispersal with monostable
nonlinearity}


\author[mpi]{J\'er\^ome Coville},
\ead{coville@mis.mpg.de}
\author[dim]{Juan D\'avila\corauthref{cor}},
\corauth[cor]{Corresponding author.} 
\ead{jdavila@dim.uchile.cl}
\author[dim]{Salom\'e Mart\'\i nez}
\ead{samartin@dim.uchile.cl}

\address[mpi]{Max-Planck Institute for Mathematics in the Sciences,\\
Inselstrasse 22\\
D-04103     Leipzig, Germany}

\address[dim]{Departamento de Ingenier\'ia Matem\'atica\\ and Centro de Modelamiento Matem\'atico(UMI CNRS 2807)\\
Universidad de Chile\\
Blanco Encalada 2120 - 5 Piso\\
Santiago - Chile}

\begin{abstract}
We study the travelling wave problem
$$
J\star u -u -cu^{\prime}+f(u)=0\quad{ in } \quad \R, \quad
u(-\infty)=0, \quad u(+\infty)=1
$$
with an asymmetric kernel $J$ and a monostable nonlinearity. We
prove the existence of a minimal speed, and under certain
hypothesis the uniqueness of the profile for $c\not=0$. For $c=0$
we show examples of non-uniqueness.
\end{abstract}

\begin{keyword}
integral equation \sep nonlocal anisotropic dispersal \sep travelling waves \sep KPP nonlinearity

\end{keyword}
\end{frontmatter}

\section{Introduction and main results}
During the past ten years, much attention has been drawn to the study of the following nonlocal equation
\begin{align}
&\frac{\partial U}{\partial t}=\j\star U -U +f(U) \quad{ in }\quad
\R^n\times \R^+ \label{mono.eq.noloc},
\\
&U(x,0)=U_0(x)
\end{align}
where $\j$ is a probability density on $\R^N$ and $f$ a given nonlinearity. Such kind of equations appear in various applications  ranging from population dynamics to Ising models as seen in \cite{BFRW,CF,DGP,EMc,F2,FoMe,HMMV,M,Sch}  among many
references. Here we will only be concerned with probability densities $\j$ which satisfy the following assumption:
\begin{align*}
\j\in C(\R^n),\, \j(z)\ge 0,\, \int_{\R^n} \j(z)
dz=1,\, \int_{\R^n} |z| \j(z)\,dz<\infty,
\end{align*}
and nonlinearities $f$  of monostable type, e.g.
\begin{align}
\tag{f1} \label{f1}
\begin{aligned}
& f\in C^1(\R),\text{ which satisfies } f(0)=f(1)=0,\,
f'(1)<0,\,f|_{(0,1)}>0
\\
& \text{ and } f|_{\R\setminus [0,1]}\le 0.
\end{aligned}
\end{align}
Such nonlinearities are commonly used in population dynamics to
describe the interaction (birth, death , \ldots ) of a species in
its environment as described in \cite{F1,GK}.


Our analysis in this paper will mainly focus on the travelling
wave solutions of  equation \eqref{mono.eq.noloc}
These particular type of solutions  are
of the form $U_e(x,t):=u(x.e+ct)$ where $e\in \mathbb{S}^{n-1}$ is
a given unit vector, the velocity $ c\in\R$ and the scalar
function $u$ satisfy
\begin{align}
&J\star u -u -cu^{\prime}+f(u)=0\quad{ in } \quad
\R,\label{mono.eq.tw}
\\
&u(-\infty)=0,\label{mono.eq.twbc-}\\
&u(+\infty)=1,\label{mono.eq.twbc+}
\end{align}
where  $u(\pm\infty)$ denotes the limit of $u(x)$ as
$x\to\pm\infty$ and $J$ is the real function defined as
$$
J(s):=\int_{ \Pi_s } \j(y)\,dy 
$$
where $ \Pi_s = \{ y \in \R^N : \langle y,e \rangle=s \}$.
Thus we shall assume that the kernel $J$ satisfies
\begin{align}
\tag{j1} \label{j1}  J\in C(\R),\, J(z)\ge 0,\, \int_{\R} J(z)
dz=1,\, \int_{\R} |z| J(z)\,dz<\infty.
\end{align}
We will call a solution $u \in L^\infty(\R)$ to
\eqref{mono.eq.tw}--\eqref{mono.eq.twbc+} a travelling wave or
travelling front if it is non-decreasing.

The first  works  to study  travelling fronts in this setting are
due  to  Schumacher \cite{Sch} and in related nonlocal problems by Weinberger \cite{W1,W2} who
constructed  travelling fronts satisfying some exponential decay
for $J$ symmetric and particular monostable nonlinearities, the so
called KPP nonlinearity, e.g.
\begin{align}
\tag{f2} \label{f2}  \text{ $f$ is monostable and satisfies
$f(s)\le f'(0)s$.}
\end{align}
Then, Harris, Hudson and Zinner \cite{HHZ} and more recently Carr
and Chmaj \cite{CC} ,  Chen and Guo \cite{CG} and Coville and
Dupaigne \cite{CD2} extended and completed the work of Schumacher
to  more general monostable nonlinearities and dispersal kernels
$J$ satisfying what is called in the literature the Mollison
condition \cite{KM,Mo,M}:
\begin{align}
\tag{j2} \label{j2} \exists \lambda>0 \quad\text{  such that }\quad \int_{-\infty}^{\infty}J(-z)e^{\lambda z}\,dz <+\infty.
\end{align}
More precisely,  they show that

\medskip
\begin{theorem} \label{thm symmetric trav wave}
\cite{CC,CG,CD2,HHZ,Sch}
Let $f$ be a monostable nonlinearity,  $J$ be a symmetric function
satisfying \eqref{j1}-\eqref{j2}. Then there exists a constant
$c^*>0$ such that for all $c\ge c^*$, there exists an increasing
function $u$, such that $(u,c)$ is a solution of
\eqref{mono.eq.tw}- \eqref{mono.eq.twbc+} and for any $c<c^*$,
there exists no increasing solution of \eqref{mono.eq.tw}-
\eqref{mono.eq.twbc+}. Moreover, if in addition $f'(0)>0$, then
any bounded  solution $(u,c)$ of \eqref{mono.eq.tw}-
\eqref{mono.eq.twbc+} is unique up to translation.
\end{theorem}

\medskip

Furthermore, as in the classical case, when the nonlinearity is
KPP the critical speed $c^*$ can be precisely evaluated by means
of a formula.

\medskip
\begin{theorem} \label{thm symmetric speed}
\cite{CC,CG,HHZ,Sch,W2} Let $f$ be a KPP nonlinearity and $J$ be
a symmetric function satisfying \eqref{j1}-\eqref{j2}. Then the
critical speed $c^*$ is given by
\begin{align*}
c^* = \min_{\lambda>0} \frac{1}{\lambda} \left( \int_\R J(x)
e^{\lambda x} \, d x + f'(0) -1 \right).
\end{align*}
\end{theorem}

In Theorems~\ref{thm symmetric trav wave} and \ref{thm symmetric
speed} the dispersal kernel $J$ is assumed to be symmetric. This
corresponds to the situation where the dispersion of the species
is isotropic.  Since the dispersal of an individual can be
influenced in many ways (wind, landscape,\ldots), it is natural to
ask what happens when the kernel $J$ is non-symmetric.
In this direction, we have the following result:


\medskip
\begin{theorem}\label{mono.th.sup}
Let $f$ be a monostable nonlinearity satisfying \eqref{f1} and $J$
be a dispersal kernel satisfying \eqref{j1}. Assume further that
there exists $(w,\kappa)$ with $w\in C(\R)$ a super-solution of
\eqref{mono.eq.tw}-\eqref{mono.eq.twbc+} in the sense:
\begin{equation}
\label{supersolution}
\begin{array}{l}
J\star w -w -\kappa w^{\prime}+f(w)  \le 0\quad{ in } \quad \R,\\
w(-\infty)\ge 0,\\
w(+\infty)\ge 1
\end{array}
\end{equation}
and such that  $w (x_0)<1$ for some $x_0 \in \R$. Then there
exists a critical speed $c^*\le\kappa$, such that for all $c\ge
c^*$ there exists a non decreasing solution $(u,c)$ to
\eqref{mono.eq.tw}-\eqref{mono.eq.twbc+} and for $c<c^*$there
exists no non-decreasing travelling wave with speed $c$.
\end{theorem}

We emphasize that in the above theorem  we do not require
monotonicity of  the supersolution $w$. The first  consequence of
Theorem \ref{mono.th.sup} is to relate  the existence of a minimal
speed $c^*$ and the existence of a travelling front for any speed
$c \ge c^*$ to the existence of a supersolution. In other words,
we have  the following  necessary and sufficient condition:

\medskip
\begin{corollary}
Let $f$ and $J$  be such that \eqref{f1} and \eqref{j1} holds.
Then there exists a non decreasing solution with minimal speed
$(u,c^*)$ of \eqref{mono.eq.tw}-\eqref{mono.eq.twbc+}    if and
only if there exists a supersolution $(w,\kappa)$ of
\eqref{mono.eq.tw}-\eqref{mono.eq.twbc+}.
\end{corollary}
\medskip

The existence of a supersolution in Theorem~\ref{mono.th.sup} is
automatic under extra assumptions on $J$. For instance, we have

\medskip


\begin{theorem}\label{mono.th.mol}
Let $f$ be a monostable nonlinearity and  $J$ satisfy \eqref{j1} and Mollison's condition \eqref{j2}. 
Then there exists a critical speed $c^*\le\kappa$, such that for all $c\ge
c^*$ there exists a non decreasing function $u$ such that $(u,c)$
is a solution of \eqref{mono.eq.tw}-\eqref{mono.eq.twbc+}. While there is no non decreasing travelling wave with speed
$c<c^*$.
\end{theorem}

\medskip

Next we examine the validity of Theorem~\ref{thm symmetric speed}
for nonsymmetric $J$. Let $c^1$ denote the following quantity

$$
c^1:=\inf_{\lambda>0} \frac{1}{\lambda} \left( \int_\R J(-x)
e^{\lambda x} \, d x + f'(0) -1 \right) .
$$
For $c\ge c^1$ we  denote $\lambda(c)$ the unique minimal
$\lambda>0$ such that
$$
-c\lambda +\int_{\R} J(-x)e^{\lambda x}dx+f'(0)-1=0.
$$
We generalize a result of Carr and Chmaj \cite{CC} to the case
when $J$ is nonsymmetric.

\medskip
\begin{theorem} \label{thm exp decay}
Let $f$ be a monostable nonlinearity satisfying \eqref{f1},
$f'(0)>0$, $f\in C^{1,\gamma}$ near $0$ and there is $m \ge1$,
$\delta>0$, $A>0$ such that
\begin{align}\label{holder}
|u-f(u)|\ge A u^m \quad \hbox{ for all $0\le u<\delta$.}
\end{align}
Let $J$ be a dispersal kernel satisfying \eqref{j1}, $J\in C^1$
and is compactly supported. Then $c^1 \le c^*$. Moreover, if $u$
is a solution of \eqref{mono.eq.tw}, \eqref{mono.eq.twbc-},  $0\le
u\le 1$, $u \not \equiv 0$ then, when $c=c^1$
\begin{equation}\label{behavior1}
0<\lim_{x\to-\infty} \frac{u(x)}{ |x| e^{\lambda (c^*) x}}
<\infty,
\end{equation}
and when $c>c^1$
\begin{equation}\label{behavior2}
0<\lim_{x\to-\infty} \frac{u(x)}{ e^{\lambda (c) x}}  <\infty.
\end{equation}
\end{theorem}

In Theorem~\ref{thm exp decay} we do not need to assume that the
solution $u$ to \eqref{mono.eq.tw}, \eqref{mono.eq.twbc-} is
monotone.

\medskip
\begin{corollary} \label{corollary velocity kpp}
If $f$ and $J$ satisfy the hypotheses of Theorem~\ref{thm exp
decay} and $f$ satisfies also \eqref{f2} then
$$
c^* = c^1.
$$
\end{corollary}

Observe that when $J$ is symmetric, by Jensen's inequality $c^1>0$. On the other hand, it is not difficult to construct examples of nonsymmetric $J$ such that $c^1 \le 0$. 
This fact should not be surprising. Indeed,  let us recall a connection between the nonlocal problem \eqref{mono.eq.noloc} and a local version which arises by considering a family of kernels that approaches a Dirac mass, that is, $J_\eps (x)= \frac{1}{\eps} J(\frac{x}{\eps})$ with $\eps>0$.  Assuming that $u$ is smooth and $J$ decays fast enough, expanding $J_\eps \star u - u $ in powers of $\eps$ we see that
\begin{align}
\nonumber
J_\eps\star u (x)- u (x) &= \frac{1}{\eps} \int_\R J(\frac{x-y}{\eps}) ( u(y) - u(x) ) \, d y =
\int_\R J(-z) ( u ( x + \eps z ) - u(x) ) \, d z 
\\
\label{expansion}
& = \eps \beta u'(x) + \eps^2 \alpha u''(x) + o(\eps^2)
\end{align}
as $\eps\to 0$, where
$$
\alpha = \frac{1}{2} \int_\R J(z) z^2 \, d z \quad \hbox{and} \quad
\beta = \int_\R J(-z) z \, d z.
$$
Thus there is a formal analogy between $J\star u - u $ and $  \beta u'(x) + \eps \alpha u''(x) $. When $J$ is symmetric then $\beta=0$ and the results for travelling waves of \eqref{mono.eq.tw}--\eqref{mono.eq.twbc+} are similar to those for travelling wave solutions of 
\begin{align}
\label{local trav wave}
\tilde \alpha u'' - c  u' + f(u) &= 0 \quad \hbox{in $\R$}, \quad u(-\infty) = 0, \quad u(+\infty)=1 ,
\end{align}
where $\tilde \alpha>0$. For \eqref{local trav wave} there exists a minimal speed $c^*>0$ such that travelling front solutions exist if and only if $c \ge c^*$ (see \cite{KPP}). For general asymmetric $J$ we see  from \eqref{expansion} that a better analogue than \eqref{local trav wave} for \eqref{mono.eq.tw}--\eqref{mono.eq.twbc+} is the problem
$$
\tilde \alpha u'' - (c- \tilde \beta) u'  + f(u) =0\quad \hbox{in $\R$}, \quad u(-\infty)=0, \quad u(+\infty)=1
$$
for some $\tilde \alpha \ge 0$ and $\tilde \beta \in \R$. This equation is the same as \eqref{local trav wave} with a shift in the speed, that is, the minimal speed is $c^* + \tilde \beta$ where $c^*$ is the old minimal speed in \eqref{local trav wave}. This new minimal speed can be either positive or negative depending on the size and sign of $\tilde \beta$, which is related to the asymmetry of $J$.

\medskip
Regarding the uniqueness of the profile of the travelling waves we
prove:

\medskip
\begin{theorem}\label{mono.th.uniq}
Assume $f$ and $J$ satisfy the hypotheses of Theorem~\ref{thm exp
decay} and $J$ satisfies:
\begin{align}
\label{MP} \exists\; a<0<b  \quad \text{ such that } \quad
J(a)>0, \; J(b)>0.
\end{align}
Then for $c\ne 0$ the   solution of the problem
\eqref{mono.eq.tw}-\eqref{mono.eq.twbc+}  is unique up to
translation.
\end{theorem}

\medskip

We notice if $c\not=0$ then any solution to \eqref{mono.eq.tw} is
continuous. In the case $c=0$, the same argument used to prove
Theorem~\ref{mono.th.uniq} gives uniqueness for continuous
solutions of \eqref{mono.eq.tw}-\eqref{mono.eq.twbc+} provided
that this problem admits a continuous solution (see
Remark~\ref{one continuous}). In the case $c=0$ one sufficient
condition for a solution $0\le u \le 1$ to \eqref{mono.eq.tw} to
be continuous is that
$$
u-f(u) \hbox{ is strictly increasing in $[0,1]$.}
$$
In Proposition~\ref{not unique} we give examples of $f$ and
non-symmetric $J$ such that no solution of
\eqref{mono.eq.tw}-\eqref{mono.eq.twbc+} is continuous, and this
problem admits infinitely many solutions.

Our results also have implications in the study of solutions to
\begin{equation} \label{nonlocal ind of x}
J\star u -u + f(u)=0
\end{equation}
which corresponds to \eqref{mono.eq.tw} with velocity $c=0$. In
\cite{CDM} it was shown that if $f(u)/u$ is decreasing and $J$ is
symmetric then any non-trivial bounded solution of \eqref{nonlocal
ind of x} is identically 1. The symmetry of $J$ was important in
the argument and it was conjectured that if the kernel $J$ is not
even \equ{nonlocal ind of x} may have more than one solution. For
this discussion we shall assume that $f$ and $J$ satisfy the
hypotheses of Theorem~\ref{thm exp decay} and $f$ also satisfies
\eqref{f2}. We observe that when the dispersal kernel is not even,
the critical velocity $c^*$ can be non-positive. If $c^*\le 0$ we
obtain that the equation \eqref{nonlocal ind of x} has a
non-constant positive  solution satisfying
\equ{mono.eq.twbc-}-\equ{mono.eq.twbc+}. Similarly, equation
\equ{nonlocal ind of x} has positive solutions satisfying
$$
\lim_{x\to-\infty} u(x) = 1 , \quad \lim_{x\to+\infty} u(x) = 0,
\quad  \hbox{$u$ is non-increasing}$$ if and only if $c_{*} \le 0$
where
\begin{equation*} 
c_{*}= \min_{\lambda>0} \frac{1}{\lambda} \left( \int_\R J(x)
e^{\lambda x} \, d x + f'(0) -1 \right).
\end{equation*}
Observe that by Jensen's inequality we have $c^* >0 $ or $c_*>0$. In summary, besides $u\equiv0$ and $u\equiv 1$ equation \eqref{nonlocal ind of x} has travelling wave solutions if $c^* \le 0$ or $c_* \le 0$. One may wonder whether other types of solutions may exist, maybe not monotone or with other behavior at $\pm\infty$. Under some additional conditions on $f$ we have a complete classification result for \eqref{nonlocal ind of x}, in the sense that we do not require the boundary conditions at $\pm \infty$, continuity nor the monotonicity of the solutions. This result can be shown by slightly modifying the arguments for Theorem~2.1 in \cite{CC}.


\medskip

\begin{theorem}\label{theo-uniquenessCC} Suppose $f$ and $J$
satisfy the hypotheses of Theorem~\ref{thm exp decay}, $J$
satisfies \eqref{MP} and  $f'(r) \le f'(0)$ for  $r \in ( 0,1)$.
Then any solution $0\le u \le 1$ of problem \eqref{nonlocal ind of
x} is one of the following: 1) $u\equiv 0$ or $u\equiv 1$, 2) a
non-decreasing travelling wave or 3) a non-increasing travelling
wave. Moreover in cases 2) and 3) the profile is unique up to
translation.
\end{theorem}

\medskip

Regarding Mollison's condition \eqref{j2} let us mention that
recently Kot and  Medlock in \cite{KM} have  shown that for a one
dimensional problem  when the dispersal kernel $J$ is even with a
fat tails and $f(s):=s(1-s)$,  the solutions of the initial value
problem \eqref{mono.eq.noloc} do not behave like travelling waves
with constant speed but rather like what they called {\it
accelerating} waves. Moreover, they predict the apparition of
accelerating waves for \eqref{mono.eq.noloc}. More precisely,
supported by numerical evidence and analytical proof, they
conjecture that \eqref{mono.eq.noloc} admits travelling wave
solutions if and only if for some $\lambda>0$
$$
\int_{-\infty}^{+\infty} J(z)e^{\lambda z}\,dz<+\infty.
$$
It appears from our analysis on non symmetric dispersal kernels,
that the existence of travelling waves with constant speed is more
related to
\begin{equation*} 
\int_{0}^{+\infty} J(z)e^{\lambda z}\,dz<+\infty \qquad\text{ for
some }\; \lambda>0
\end{equation*}
if we look at fronts propagating from the left to the right and
\begin{equation*} 
\int_{0}^{+\infty} J(-z)e^{\lambda z}\,dz<+\infty \qquad\text{ for
some }\; \lambda>0
\end{equation*}
if we look at  fronts propagating from the right to the left. As a
consequence, for asymmetric kernels, it may happen that in one
direction, the solution behave like a front with finite speed and
in the other like an accelerating wave.

\medskip

The outline of this paper is the following. In Section
\ref{mono.s.pre}, we recall some results on front solutions for
ignition nonlinearities,  then in Section \ref{mono.s.constrsi} we
construct increasing solution of  for $J$ compactily supported.     Section \ref{mono.s.pthsup} is devoted
to the proofs of Theorem \ref{mono.th.sup} and 
Theorem \ref{mono.th.mol}. Section~\ref{velocity}
contains the proof of Theorem~\ref{thm exp decay} and
Corollary~\ref{corollary velocity kpp}. In Section~\ref{section
uniqueness} we prove the uniqueness of the profile
Theorem~\ref{mono.th.uniq} and Theorem~\ref{theo-uniquenessCC}.


\section{Approximation by ignition type nonlinearities}
\label{mono.s.pre}

The proof of Theorem \ref{mono.th.sup} essentially relies on some
estimates and properties of the speed of fronts for problem
\eqref{mono.eq.noloc} with ignition type nonlinearities $f$. We
say that $f$ is of ignition type if $f\in C^{1}([0,1])$ and
\begin{align}
\tag{f3}\label{f3}
\begin{aligned}
&\hbox{ there exists $\rho\in(0,1)$ such that $f_{|[0,\rho]}\equiv
0$,$f_{|(\rho,1)}>0$ and $f(1)=0$.}
\end{aligned}
\end{align}

Consider the following problem
\begin{align}\label{mono.eq.pre epsilon 0}
\left \{
\begin{array}{rcl}
J\star u - u -cu' +f(u)&=&0  \quad \text{ in $\R$} \\
u(-\infty)&=&0  \\
u(+\infty)&=&1,
\end{array}
\right.
\end{align}
where $c\in\R$ and  $f$ is either  an ignition nonlinearity or a
monostable nonlinearity.

The main result in this section is the following:
\begin{proposition}\label{mono.th.m}
Let $f$ be a monostable nonlinearity and assume that $J$ is a non
negative continuous function of  unit mass. Assume further that
there exists $(w,\kappa)$  a super-solution of
\eqref{mono.eq.tw}-\eqref{mono.eq.twbc+}. Let $(f_k)_{k\in\N}$ be
any sequence of ignition functions which converges pointwise to
$f$ and satisfies $\forall k\in\N, f_k\le f_{k+1} \le f$ and let
$c_k$ be the unique speed of fronts associated to
\eqref{mono.eq.pre epsilon 0}. Then
\begin{equation}\label{mono.eq.clim}
\lim_{k\to+\infty}c_k=c^*,
\end{equation}
exists and is independent of the sequence $f_k$. Furthermore,
$c^*\le \kappa$, there exists a non decreasing solution $(u,c^*)$
of \eqref{mono.eq.tw}-\eqref{mono.eq.twbc+} and for $c<c^*$ there
are no non-decreasing solutions to
\eqref{mono.eq.tw}-\eqref{mono.eq.twbc+}.
\end{proposition}

The fact that for \eqref{mono.eq.pre epsilon 0} with ignition type
nonlinearity there exists a unique speed of fronts has been
recently established by one of the authors in \cite{Co2,Co3,Co4}
and holds also for the following perturbation of
\eqref{mono.eq.pre epsilon 0}
\begin{align}\label{mono.eq.pre}
\left \{
\begin{array}{rcl}
\epsilon u'' + J\star u - u -cu' +f(u)&=&0  \quad \text{ in $\R$} \\
u(-\infty)&=&0  \\
u(+\infty)&=&1,
\end{array}
\right.
\end{align}
where  $\epsilon\ge 0$, $c\in\R$.

\begin{theorem} (\cite[Theorem 1.2]{Co4} and \cite[Theorem 3.2]{Co2})
\label{mono.th.ig} Let $f$ be an ignition nonlinearity and assume
that $J$ satisfies \eqref{j1}.  Then there exists a non decreasing
solution $(u,c)$ of \eqref{mono.eq.pre}. Furthermore the speed $c$
is unique. Moreover, if $(v,c')$ is a super solution of
\eqref{mono.eq.pre}, then $c\le c'$. The inequality becomes strict
when $v$ is not a solution of \eqref{mono.eq.pre}.
\end{theorem}
We remark that in this results the supersolution $v$ is not
required to be monotone.

\begin{corollary}\label{mono.cor.ig}
Let $f_1\ge f_2$, $f_1 \not\equiv f_2$ be two ignition
nonlinearities and assume that $J$ is a non negative continuous
function of  unit mass with finite first moment.  Then $c_1> c_2$
where $c_1$ and $c_2$ are the corresponding  unique speeds given
by Theorem \ref{mono.th.ig}.
\end{corollary}

We also recall some useful results  on solutions of
\eqref{mono.eq.pre}, which can be found in \cite{Co4,CD2}.
\begin{lemma} \cite[Lemma 2.1]{Co4}  \\
\label{lema 2.4} 
Suppose $f$ satisfies \eqref{f1} and $J$ satisfies
\eqref{j1}. Assume $\eps \ge 0$, $c\in \R$ and let $0 \le u \le 1$
be an increasing solution of \eqref{mono.eq.pre}. Then
$$
f(l^{\pm})=0,
$$
where $l^{\pm}$ are the limits of $u$ at $\pm\infty$.
\label{mono.lem.esti}
\end{lemma}

\begin{lemma}\cite[Lemma 2.2]{Co4} \label{mono.lem.ig}
Let $f$ and $J$ be as in Theorem~\ref{mono.th.ig}.  Then
following holds
 $$\mu c^2-\nu|c|\le 0$$
where the constants $\mu, \nu$  are defined by
\begin{align*}
&\mu:=\inf\{\rho,1-\rho\} &\nu:=\int_{\R}J(z)|z|\,dz
\end{align*}
\end{lemma}

\dem{Proof of Proposition~\ref{mono.th.m}.} Let   $(f_n)_{n\in\N}$
be a sequence of ignition functions which converges pointwise to
$f$ and satisfies $\forall n\in\N, f_n\le f_{n+1}\le f$.  Let
$(u_n,c_n)$ denote the corresponding solution given by Theorem
\ref{mono.th.ig}. By Corollary~\ref{mono.cor.ig} $(c_n)_{n\in\N}$
is an increasing sequence. Next, we see that $c_n \le \kappa$.
Since $w$ satisfies
$$
J \star w - w - \kappa w' + f_n(w) \le 0 \quad \hbox{in $\R$}
$$
by Theorem \ref{mono.th.ig} we get
$$
c_n\le \kappa .
$$

Let us observe that we can normalize the sequence of solutions
$u_n$ by $u_n(0)=\frac{1}{2}$. Indeed, when $c^*=0$ since $c_n <
c^*$ the solution $u_n$ is smooth.  Since any translation of $u_n$
is a solution of the problem and $u_n(-\infty)=0$,
$u_n(+\infty)=1$ we can normalize it by $u_n(0)=\frac{1}{2}$. When
$c^*\neq 0$, since $c_n \to c^*$ the sequence $u_n$ is smooth for
all $n$ sufficiently large. Thus the same normalization can be
also taken in this situation.

Since $(u_n)_{n\in\N}$ is an uniformly bounded sequence of
increasing functions, using Helly's lemma  there exists a
subsequence which converges pointwise to a non decreasing function
$u$. Moreover, $u$ satisfies  in the distribution sense
\begin{align*}
&J\star u - u -c^*u^{\prime} +f(u)= 0 \quad \hbox{in $\R$},
\end{align*}
and by the monotonicity and the normalization of $u_n$
\begin{align}
\label{normalization u} u(x) \le \frac{1}{2} \quad \hbox{for all
$x\le 0$}, \quad u(x) \ge \frac{1}{2} \quad \hbox{for all $x\ge
0$}.
\end{align}
Observe that  when $c^*\neq 0$, using $C^1_{loc}$ regularity, we
get that $u\in C^1_{loc}$ and satisfies the above equation in  a
strong sense. Otherwise, when $c^*=0$, a standard argument shows
that $u$ satisfies almost everywhere the  equation
\begin{align*}
&J\star u - u +f(u)= 0.
\end{align*}

Observe that by \eqref{normalization u} $u$ is non trivial. It
remains to show that $u$ satisfies the right boundary conditions.
Now, since $u$ is non decreasing and bounded, the following limits
are well defined:
\begin{align*}
&l^-:=\lim_{x\to-\infty}u(x)\\
&l^+:=\lim_{x\to+\infty}u(x).
\end{align*}
We get $l^+=1$ and $l^-=0$ using Lemma \ref{mono.lem.esti}, the
definition of $f$ and the monotonicity of $u$.

To finish we need to prove that $c^*$ is independent of the
sequence $f_n$. So consider another sequence $\tilde f_n$ of
ignition functions such that $\tilde f_n \le \tilde f_{n+1} \le f$
and $\tilde f_n \to f$ pointwise. Let $(\tilde u_n, \tilde c_n)$
denote the front solution and speed of \eqref{mono.eq.pre epsilon
0} with nonlinearity $\tilde f_n$ and let
$$
\tilde c = \lim_{n\to\infty} \tilde c_n.
$$
Since $u = \lim_{n\to \infty} u_n$ satisfies
$$
J\star u - u - c^* u' + \tilde f_n(u) \le 0
$$
by Theorem~\ref{mono.th.ig} we have $\tilde c_n \le c^*$. Hence
$\tilde c \le c^*$ and reversing the roles of $f_n$ and $\tilde
f_n$ we get $c^* \le \tilde c$.

Finally observe that for $c<c^*$ there is no monotone solution to
\eqref{mono.eq.twbc-}-\eqref{mono.eq.twbc+}. Otherwise this
solution would be a supersolution of \eqref{mono.eq.pre epsilon 0}
with $f_n$ instead of $f$. By Theorem \ref{mono.th.ig} we would
have $c_n \le c$ for all $n$, which is a contradiction.

\fdem


\section{Construction of solutions of \eqref{mono.eq.tw}-\eqref{mono.eq.twbc+} when $J$ is compactly supported}
\label{mono.s.constrsi}
In this section we construct  monotone solutions of \eqref{mono.eq.tw}-\eqref{mono.eq.twbc+} when $J$ is compactly supported. More precisely we prove the following

\medskip
\begin{proposition}\label{mono.th.jc}
Let $f$ be a monostable nonlinearity and  $J$ be   continuous
compactly supported which  satisfies \eqref{j1}. Assume further
that there exists $a\in \R$ such that $\{a,-a\}\subset supp(J)$.
Then there exists a critical speed $c^*$, such that for all $c\ge
c^*$ there exists a non decreasing function $u$ such that $(u,c)$
is a solution of \eqref{mono.eq.tw}-\eqref{mono.eq.twbc+}.
Moreover, there is no non decreasing travelling wave with speed
$c<c^*$.
\end{proposition}
\medskip

To prove the above result we proceed following the strategy developed in
\cite{CD2}. It is based on the vanishing viscosity technique,
apriori estimates,   the construction of adequate super and
sub-solutions and the characterization of the critical speed
obtained in Section \ref{mono.s.pre}. Let us first briefly explain
how we proceed.


\medskip  \noindent \textbf{Step 1:}

For convenience, let us first rewrite  problem \eqref{mono.eq.pre}
the following way:
\begin{align} \label{mono.eq.aux2}
\left \{
\begin{array}{rcl}
\opm u +f(u)&=&0  \quad \text{ in $\R$} \\
u(-\infty)&=&0  \\
u(+\infty)&=&1,
\end{array}
\right.
\end{align}
where the operator $\m$ is defined for a given $\epsilon>0$,
$c\in\R$  by
\begin{equation}\label{defM}
\begin{array}{l}
\opm u = \M(\epsilon,c)u = \epsilon u'' + J\star u - u -cu'.\\
\end{array}
\end{equation}
For problem  \eqref{mono.eq.aux2}, for small $\eps$, we construct
a super solution which is independent of $\eps$. More precisely we
show the following
\begin{lemma} \label{lemma supersolution}
Let $J$ and $f$ be as in Theorem \ref{mono.th.jc}. Then there
exists $\eps_0$ and $(w,\kappa)$ such that $\forall \eps \le
\eps_0$, $(w,\kappa)$ is a super-solution of \eqref{mono.eq.aux2}.
\end{lemma}

\medskip  \noindent
\textbf{Step 2:} Using the above super solution and a standard
approximation scheme, for fixed $0<\eps\le\eps_{0}$,  we prove the
following
\begin{proposition}\label{mono.th.jceps}
Fix $0 < \eps\le \eps_0$ and let $J$ and $f$ be as in Theorem
\ref{mono.th.jc}. Then there exists $c^*(\eps)$ such that $\forall
c\ge c^*(\eps)$, there exists an increasing function $u_\eps$ such
that $(u_\eps,c)$ is a solution of \eqref{mono.eq.aux2}. Moreover
$c^*(\eps) \le \kappa$ where $(w,\kappa)$ is the super-solution of
Lemma~\ref{lemma supersolution}.
\end{proposition}

\medskip  \noindent
\textbf{Step 3:} We study  the singular limit $\eps \to 0$ and
prove Proposition~\ref{mono.th.jc}.

\medskip

Some of the  arguments developed in  \cite{CD2}, on which this
procedure is based, do not use the symmetry of $J$. Hence in some
cases we will skip details in our proofs, making appropriate
references to \cite{CD2}.

\medskip

We divide this section in 3 subsections, each one devoted to one
Step.


\subsection{\bf Step 1. Existence of a super-solution}
\label{mono.ss.sub&sup}

We start with the construction of a super-solution of
\eqref{mono.eq.aux2} for speeds $c\ge\bar\kappa$ for some
$\bar\kappa>0$ which is independent of $\eps $ for $0<\eps \le 1$.
\begin{lemma}{\ \ }\\ \label{mono.lem.sup}
Assume $J$ has compact support and let $\eps> 0$. There exists a
real number $\bar\kappa>0$ and an increasing function $\bar w \in
C^2(\R)$ such that, given any $c\ge\bar\kappa$ and $0<\eps \le 1$
\begin{equation*}
\left\{\begin{array}{l}
\opm{\bar w} +f(\bar w)\le 0 \ \ \text{ in } \ \ \R,\\
\bar w(-\infty)=0,\\
\bar w(+\infty)=1,\\
\end{array}\right.
\end{equation*}
where $\M=\M(\eps,c)$ is defined by \eqref{defM}. Furthermore,
$\bar w(0)=\frac{1}{2}$.
\end{lemma}

The construction of the super-solution is an adaptation of the one
proposed in \cite{CD2}. The essential difference lies in the
computation of the super-solution in a neighborhood of $-\infty$.

\dem{Proof.} As in \cite{CD2}, fix positive constants $N, \lambda,
\delta$ such that $\lambda>\delta$.

Let $\bar w\in C^2(\R)$ be a positive increasing function
satisfying

\begin{itemize}
  \item $\bar w(x)=e^{\lambda x}$ for $x\in (-\infty,-N]$,
 \item $\bar w(x)\le e^{\lambda x}$ on $\R$,
 \item $\bar w(x)=1-e^{-\delta x}$ for $x\in [N,+\infty)$,
 \item $\bar w(0)=\frac{1}{2}$.
\end{itemize}
Let $x_{0}=e^{-\lambda N}$ and $x_{1}=1-e^{-\delta N}$. We have
$0<x_{0}<x_{1}<1$.\\
We now construct a positive function $g$ defined on $(0,1)$ which
satisfies $g(\bar w)\ge f(\bar w)$. Since $f$ is smooth near 0 and
1, we have for $c$ large enough, say $c\ge \kappa_0$,
\begin{equation}
 \lambda(c-\lambda)s \ge f(s) \ \ for \ \ s \in [0,x_0]
\label{mono.def.g1}
\end{equation}
and
\begin{equation}
 \delta ( c -\delta)(1-s)\ge f(s) \ \ for \ \ s \in [x_1,1].
\label{mono.def.g2}
\end{equation}
Therefore we can achieve $g(s)\ge f(s)$ for $s$ in [0,1], with $g$
defined by:
\begin{align} \label{mono.def.g}
g(s)=\left \{\begin{array}{ll}
        \lambda(\kappa_0-\lambda)s & for \ \ 0\le s \le x_0 \\
         l(s)   & for \ \ x_0< s < x_1 \\
        \delta(\kappa_0-\delta)(1-s) & for \ \ x_1\le s \le 1 \\
        \end{array}
  \right.
\end{align}
where $l$ is any smooth positive function greater than $f$ on
$[x_0,x_1]$ such that $g$ is of class $C^1$.

According to \eqref{mono.def.g}, for $x\le -N$ i.e. for $w\le
e^{-\lambda N}$, we have
\begin{eqnarray*}
\opm{\bar w}+ g(\bar w)&=&\eps \bar w'' + J\star \bar w -\bar w -
c \bar w^{\prime}+g(\bar w)
\\
&=&\eps\lambda^2e^{\lambda x} + J\star \bar w-e^{\lambda
x}-\lambda c  e^{\lambda x}+ \lambda(\kappa_0-\lambda)e^{\lambda
x}
\\
&\le& \eps\lambda^2e^{\lambda x} +  J\star e^{\lambda x}
-e^{\lambda x}-\lambda c e^{\lambda x}+
\lambda(\kappa_0-\lambda)e^{\lambda x}
\\
&\le &e^{\lambda x}[\int_{\R}J(-z)e^{\lambda z}dz -1-\lambda
(c-\kappa_{0}) -\lambda^{2}(1-\eps)]
\\
&\le& 0,
\end{eqnarray*}
for $c$ large enough, say $$c\ge
\kappa_1=\frac{\int_{\R}J(-z)e^{\lambda z}dz -1 +\lambda\kappa_{0}
-\lambda^{2}(1-\eps)}{\lambda}. $$

In the open set  $(x_1,+\infty)$, the computation of the
super-solution is identical to the one in \cite{CD2}. So, we  end
up with
$$
\opm{\bar w}+g(\bar w)\le 0 \quad \text{ in } (x_1,+\infty)
$$
for $c$ large enough, say $c\ge \kappa_2$.

Therefore,  by taking $c \ge \sup\{\kappa_0,\kappa_1,\kappa_2\}$,
we achieve
\begin{eqnarray*}
g(\bar w)\ge f(\bar w) \qquad\text{and}\qquad \opm{\bar w}+g(\bar w)\le 0\\
\text{for} \quad 0\le \bar w\le e^{-\lambda N} \quad \text{and}
\quad \bar w\ge 1-e^{-\delta N}.
\end{eqnarray*}

For the remaining values of $\bar w$, i.e. for $ x\in [-N,N]$,
$\bar w^{\prime}>0$ and we may increase $c$ further if necessary,
to achieve

$$
\opm{\bar w}+g(\bar w) \le 0 \ \ in \ \ \R.
$$

The result follows for $$\bar
\kappa(\eps):=\sup\{\kappa_0,\kappa_1,\kappa_2,\kappa_3\},$$ where
$$\kappa_3 =\sup_{x\in[-N,N]}\{\frac{\eps|\bar w''|+|J\star \bar w
-\bar w|+g(\bar w)}{\bar w'}\}.$$
 \fdem

Now, note  that $\bar\kappa(\eps)$ is a non-decreasing function of
$\eps$, therefore for all non-negative $\eps\le 1$, $(\bar w,
\bar\kappa)$ with $\bar \kappa = \bar \kappa(1)$, will be a super
solution of \eqref{mono.eq.aux2}, which ends the Step 1.

\begin{remark}
The above construction of a super-solution also works if we only
assume that for some positive $\lambda$, the following holds
$$\int_{0}^{+\infty}J(-z)e^{\lambda z}\,dz<+\infty.$$
\end{remark}


\subsection{\bf Step 2. Construction of a solution when $\eps>0$}

To prove Proposition~\ref{mono.th.jceps} we follow the strategy
used in \cite{CD2} relying on the following approximation scheme.

We first prove existence and uniqueness of a monotone solution
for 
\begin{equation}
\label{mono.eq.aux1} 
\left
\{\begin{array}{rcl}
\ops u +f(u)&=& -h_r(x) \quad \text{ in $\omega$}, \\
      u(-r)&=&\theta, \\
      u(+\infty)&=&1,   \end{array}  \right.
\end{equation}\label
where  $\epsilon>0$, $r\in\R$, $c\in\R$
and $\theta \in (0,1)$ are given, and

\begin{align}
&\omega = (-r,+\infty),\\
&\ops u = \S(\epsilon,r,c)[u] = \epsilon
u^{\prime\prime}+\int_{-r}^{+\infty}J(x-y)u(y)dy - u -
cu^{\prime},\label{defS}\\
& h_r(x)=\theta\int_{-\infty}^{-r}J(x-y)dy.
\end{align}
More precisely, we show
\begin{proposition}\label{mono.th.jcsi}
Assume $f$ and $J$ are as in Proposition~\ref{mono.th.jc}. For any
$\eps>0$, $\theta \in [0,1)$ $r>0$ so that $supp J\subset
(-r,+\infty)$ and $c\in \R$  there exists a unique positive
increasing solution $u_c$ of \eqref{mono.eq.aux1}
\end{proposition}

To prove this proposition we use a construction introduced by one
of the authors \cite{Co3,Co4} which consists first to obtain a
solution of the following problem:
\begin{equation}
  \left \{\begin{array}{l}
\opl u +f(u) +h_{r}+h_{R}= 0 \; \text{ for } x\in \O\\
u(-r)= \theta,\\
u(+R)= 1,
   \end{array}
 \right.
\label{mono.eq.aux3}
\end{equation}
where $\O=(-r,+R)$ and $\L=\L(\eps,J,r,R,c),  h_{r}$ and $h_R$ are
defined by
\begin{equation}\label{defL}
\begin{array}{l}
\opl u = \L(\epsilon,J,r,R,c)[u] = \epsilon u^{\prime
\prime}+\left[\int_{-r}^{+R}J(x-y)u(y)dy - u\right] -
cu^{\prime},\\
h_r(x)=\theta\int_{-\infty}^{-r}J(x-y)dy.\\
h_R(x)=\int_{+R}^{+\infty}J(x-y)dy.
\end{array}
\end{equation}

Namely, we have,
\begin{proposition}{\ \ }\\
Assume $f$ and $J$ are as in Proposition~\ref{mono.th.jc}. For any
$\eps>0$, $\theta \in [0,1)$ $r<R$ so that $supp J\subset (-r,R)$
and $c\in \R$  there exists a unique positive increasing solution
$u_c$ of \eqref{mono.eq.aux3}. \label{mono.th.const2}
\end{proposition}

\dem{Proof.} The construction of a solution uses the super- and
sub-solution iterative scheme  presented in \cite{Co4}. To produce
a solution, we just have to construct ordered sub and
super-solutions. An easy computation  shows that  $\under u=
\theta$ and $\bar u=1$ are respectively a sub and a super-solution
of \eqref{mono.eq.aux3}. Indeed,
\begin{align*}
\opl{\underline{u}}+f(\underline{u})+h_{r}+h_{R}
&=\int_{-r}^{R}J(x-y)\theta\,dy -\theta +
\theta\int_{-\infty}^{-r}J(x-y)\,dy
\\
& \quad  +\int_{R}^{+\infty}J(x-y)\,dy+f(\theta)
\\
&=(1-\theta)\int_{R}^{+\infty}J(x-y)\,dy +f(\theta) \ge 0
\end{align*}
and
\begin{align*}
\opl{\bar u}+f\bar u)+h_{r}+h_{R}&=\int_{-r}^{R}J(x-y)\,dy -1 +
\theta\int_{-\infty}^{-r}J(x-y)\,dy
\\
&  \quad  +\int_{R}^{+\infty}J(x-y)\,dy+f(1)
\\
&=(\theta-1)\int_{-\infty}^{-r}J(x-y)\,dy \le 0
\end{align*}
The uniqueness and the monotonicity of such solutions have been
already established in \cite{Co3}, so we refer to this reference
for interested reader. \fdem

We are now in a position to prove Proposition~\ref{mono.th.jcsi}

\dem{Proof of Proposition~\ref{mono.th.jcsi}.} Let us now
construct a solution  of \eqref{mono.eq.aux1}. Fix $\eps>0$, $c\in
\R$ and $r>0$ such that $supp(J)\subset \o$. Let $(R_n)_{n\in \N}$
be a sequence of real which converges to $+\infty$. Since $J$ has
compact support, without loosing generality we may also assume
that $supp(J)\subset (-r,R_n)$ for all $n\in \N$. Let us denote
$(u_n,c)$ the corresponding solution given by Proposition
\ref{mono.th.const2}. Clearly, $h_{R_n}\to 0$ pointwise, as
$n\to\infty$. Observe now that $(u_n)_{n\in\N}$ is a uniformly
bounded sequence of increasing functions. Since $\eps>0$, using
local $C^{2,\alpha}$ estimates, up to a subsequence, $u_n$
converges in $C^{2,\alpha}_{loc}$ to a non-decreasing function
$u$. Therefore $u\in C^{2,\alpha}$ and  satisfies 
\begin{equation}
\label{mono.eq.sol} 
\left\{\begin{array}{ll} \epsilon u^{\prime
\prime}+\int_{-r}^{+\infty}J(x-y)u(y)\,dy - u - cu^{\prime}
+f(u)+h_r= 0 \ \ \text{ in } \ \ \o
\\
u(-r)= \theta
\end{array}
\right. 
\end{equation}
To complete the construction of the
solution, we prove that $u(+\infty)= 1$. Indeed, since $u$ is
uniformly bounded and non-decreasing, $u$ achieves its limit at
$+\infty$. Using Lemma \ref{mono.lem.esti} yields $u(+\infty)= 1$.
\fdem

\dem{Proof of Proposition~\ref{mono.th.jceps}.} By
Lemma~\ref{mono.lem.sup} there exists $\bar \kappa$ and a function
$\bar w$ which is a supersolution to \eqref{mono.eq.aux2} for any
$c \ge \bar \kappa$ and any $0<\eps\le 1$. If $c \ge \bar \kappa$,
following the approach in \cite{CD2}, we can take the limit as
$r\to \infty$ in the problem \eqref{mono.eq.aux1} to obtain a
solution of \eqref{mono.eq.aux2}.

Finally one can also verify, see \cite{CD2}, that there exists a
monotone solution $u_\eps$ with the following speed
$$
c^*(\eps):=\inf\{\, c\,| \, \eqref{mono.eq.aux2} \; \text{ admits
a monotone solution with speed }\; c\}.
$$
The proof of these claims are straightforward adaptations of
\cite{CD2}, since in this reference the author makes no use of the
symmetry of $J$ for this part of the proof, and essentially relies
on the Maximum principle and Helly's Theorem. We point the
interested reader to \cite{CD2} for the details. \fdem

\begin{remark}\label{mono.rm.esticeps}
Note that from the previous comments we get the following uniform
estimates
$$
\forall \, 0<\eps\le\eps_0\quad c^*(\eps)\le \bar \kappa .
$$
\end{remark}


\subsection{\bf Step 3. Proof of Proposition~\ref{mono.th.jc} \label{mono.ss.step3}}

We essentially use the ideas introduced in \cite{CD2}.

First, we remark that since  $J$ has a compact support, using the
super-solution of Step~1, we  get  from
Proposition~\ref{mono.th.m} a monotone solution $(u,c^*)$ of
\eqref{mono.eq.tw} -- \eqref{mono.eq.twbc+}. Furthermore, there
exists no monotone solution of  \eqref{mono.eq.tw} --
\eqref{mono.eq.twbc+} with speed $c<c^*$ and we have the following
characterization:
$$
\lim_{k \to \infty} c_{k}=c^*,
$$
where $c_k$ is the unique speed of fronts associated with an
arbitrary sequence of ignition functions $(f_k)_{k\in\N}$ which
converges pointwise to $f$ and satisfies $\forall k\in\N, f_k\le
f_{k+1} \le f$.

Also observe that from Remark \ref{mono.rm.esticeps} we have a
uniform bound from above on $c^*(\eps)$.
\begin{lemma}
For all $\eps \le \eps_0$ we have $c^*(\eps)\le \bar \kappa$.
\end{lemma}
For any speed $c\ge \bar \kappa > 0$, there exists a monotone
solution $(u_\eps,c)$ of \eqref{mono.eq.aux2} for any $\eps\le
\eps_0$. Normalizing the functions by $u_{\eps}(0)=\frac{1}{2}$
and letting $\eps \to 0$, using Helly's Theorem, a priori bounds
and some regularity we end up with a solution $(u,c)$ of
\eqref{mono.eq.tw} -  \eqref{mono.eq.twbc+}. Repeating this
limiting process for any speed $c\ge \bar \kappa$, we end up with
a monotone solution of \eqref{mono.eq.tw} - \eqref{mono.eq.twbc+}
for any speed $c\ge \bar \kappa$.

Define now the following critical speed
$$
c^{**}=\inf \{\, c\, |\, \forall\, c^\prime \geq c \ \
\eqref{mono.eq.tw}-\eqref{mono.eq.twbc+}\ \ \text { has a positive
monotone solution of speed }\ \  c^\prime\}.
$$
\begin{remark}
Observe that from the uniform bounds we easily see that
\begin{align}
\label{eq.climinf} c^{**}\le \liminf_{\eps\to 0} c^*(\eps).
\end{align}
\label{mono.rm.climinf}
\end{remark}
Obviously, we have $c^*\le c^{**} \le \bar \kappa$. To complete
the proof of Proposition~\ref{mono.th.jc}, we are then led to prove
that $c^{**}=c^*$. To prove  this equality, we use some
properties of the speed of the following approximated problem
\begin{equation}
\label{mono.eq.aux4}
\left \{
\begin{array}{ll}
\eps u''+   J\star u - u - cu^{\prime} +f\eta_{_{\theta}}(u)= 0 &
\text{in} \  \ \R,
\\
u(-\infty)=0,\\
u(+\infty)=1,
\end{array}
\right. 
\end{equation} 
where $\theta>0$, $\eta_{_{\theta}}(u)
= \eta(u/\theta)$ and $\eta \in C^\infty(\R)$ is such that
$$
\hbox{$0 \le \eta \le 1$, \quad $\eta'\ge0 $, \quad $\eta(s) =0$
for $s\le 1$, \quad $\eta(s)=1$ for $s\ge 2$.}
$$
Then $\eta_{_{\theta}}$ has the following properties
\begin{itemize}
\item
$\eta_{_{\theta}} \in C^{\infty}(\R)$,
\item
$0\le \eta_{_{\theta}}\le 1$,
\item
$ \eta_{_{\theta}}(s)\equiv 0$ for $s\le \theta$ and  $
\eta_{_{\theta}}(s)\equiv  1$ for $s\ge 2\theta$,
\item if $0 <  \theta_1 \le \theta_2$ then
$\eta_{_{\theta_1}} \le \eta_{_{\theta_2}}$.
\end{itemize}

\medskip
For  \eqref{mono.eq.aux4}, we have the following results:
\begin{lemma}\label{limits}
Let $c^{\theta}$ be the unique speed of front solutions to
\eqref{mono.eq.pre epsilon 0} and nonlinearity $f
\eta_{_{\theta}}$. Let $c_\eps^{\theta},c^*(\eps)$ be respectively
the unique (minimal) speed solution of \eqref{mono.eq.aux2} with
the nonlinearity $f \eta_{_{\theta}}$ and $f$. Then the following
holds:
\begin{itemize}
\item[a)]
For fixed $\theta > 0$, \quad $\lim_{\eps\to
0}c_\eps^{\theta}=c^{\theta}.$
\item[b)]
For fixed $\eps$ so that $\eps_0\ge \eps>0$, \quad
$\lim_{\theta\to 0}c_\eps^{\theta}=c^*(\eps).$
\end{itemize}
\label{mono.lem.ctheta}
\end{lemma}
\dem{Proof.} The first limit, as
$\eps \to 0$ when $\theta>0$ is fixed, has been already obtained
in \cite{Co4}, so we refer to this reference for a detailed proof.
The second limit, for fixed $\eps>0$, is obtained using a similar
argument as in the proof of Proposition~\ref{mono.th.m} to obtain
the characterization of $c^*$. \fdem

\dem{Proof of Proposition~\ref{mono.th.jc}.} Assume by contradiction
that $c^*<c^{**}$. Then choose $c$ such that $ c^{*}< c< c^{**}$.
By \eqref{eq.climinf} we may fix $\eps_0>0$ small such that
\begin{align} \label{asd}
c < c^*(\eps) \quad \forall \eps \in (0,\eps_0).
\end{align}
Now consider any sequence $\bar \theta_n \to 0$. Since $c^{\bar
\theta_n}< c$, using Lemma~\ref{limits}~a) there exists
$0<\eps_n<\eps_0$, $\eps_n\to 0$ such that
\begin{align}
\label{left} c_{\eps_n}^{\bar \theta_n} < c .
\end{align}
Then, using the continuity of the map $\theta \mapsto c_{\eps_n}^
\theta$, \eqref{left}, \eqref{asd} and Lemma~\ref{limits}~b) we
conclude that there exists $0<\theta_n<\bar \theta_n$ such that
$$
c = c_{\eps_n}^{\theta_n}.
$$
Note that $\theta_n\to 0$. Let $u_n$ be the associated solution to
\eqref{mono.eq.aux2} with $\eps = \eps_n$, speed $c$ and
nonlinearity $f \eta_{_{\theta_n}}$. We normalize $u_n$ by
$u_n(0)=1/2$. Using Helly's theorem we get a solution $\bar u$ of
\eqref{mono.eq.tw}-\eqref{mono.eq.twbc+} with speed $c$. This
contradicts the definition of $c^{**}$.
\fdem


\section{Construction of solution in the general case: proof of Theorems~\ref{mono.th.sup}  and \ref{mono.th.mol}}
\label{mono.s.pthsup}

Theorem~\ref{mono.th.mol} is a direct consequence of Theorem~\ref{mono.th.sup}.  Indeed,  since $J$ satisfies the Mollison condition, the construction  in Section~\ref{mono.s.constrsi} (Step 1, subsection~\ref{mono.ss.sub&sup}) of a smooth super-solution $(w,\kappa)$ with $w(0)=\frac{1}{2}$ holds. Therefore,  Theorem~\ref{mono.th.mol} is a  direct application of Theorem~\ref{mono.th.sup}.

In the rest of the section we prove  Theorem~\ref{mono.th.sup}, that is, 
we construct solutions of \eqref{mono.eq.tw}-\eqref{mono.eq.twbc+} only
assuming that  there exists a super-solution $(w,\kappa)$ of
\eqref{mono.eq.tw}-\eqref{mono.eq.twbc+}.
The construction uses a standard procedure of
approximation of $J$ by kernels $J_{n}$ with compact support and
the characterization of the minimal speed $c^*$ obtained in
Section \ref{mono.s.pre}.

Let us describe briefly our proof. From
Proposition~\ref{mono.th.m}, there exists a monotone solution
$(u,c^*)$ of \eqref{mono.eq.tw}-\eqref{mono.eq.twbc+} with
critical speed. Then we construct monotone solution of
\eqref{mono.eq.tw}-\eqref{mono.eq.twbc+}  for any $c>c^*$, $c\neq
0$, using a sequence $(J_n)_{n\in\N}$ of  approximated kernels and
the same type of arguments developed in the Step 3 of the above
section. Let us first construct the approximated kernel and get
some uniform lower bounds for $c_n^*$.


\subsection{\bf The  approximated kernel and related problems}

First, let $j_0$  be a positive symmetric function
defined by
\begin{equation}
j_0(x)=\left\{\begin{array}{lc}
e^{\frac{1}{x^2-1}}& \text{ for } x\in (-1,1)\\
0 & \text{ elsewhere }
\end{array}
\right.
\end{equation}
Now, let $(\chi_n)_{n\in\N}$ be the following  sequence of
``cut-off'' function:
\begin{itemize}
\item
$\chi_{n}\in C_{0}^{\infty}(\R)$,
\item
$0\le \chi_{n}\le 1$,
\item
$ \chi_{n}(s)\equiv 1$ for $|s|\le n$ and  $ \chi_{n}(s)\equiv  0$
for  $|s|\ge 2n$.
\end{itemize}

Define
$$J_n:=\frac{1}{m_n}\left(\frac{j_0}{n}+J(z)\chi_{n}(z)\right),$$
where $m_n:= \frac{1}{n}\int_{\R}j_0(z)dz
+\int_{\R}J\chi_{n}(z)\,dz$. Observe that since $\int_{\R}j_0>0$,
$J_n$ is well defined and $J_n(z)\to J(z)$ pointwise.

Since $J_n$ satisfies the assumption of Proposition~\ref{mono.th.jc},
there exists for each $n\in \N$ a critical speed $c^*_n$ for the
problem \eqref{mono.eq.japprox} below:
\begin{equation}
\label{mono.eq.japprox} 
\left\{
\begin{array}{l}
J_n\star u -u -cu^{\prime} +f(u)= 0\qquad\text{  in } \R
\\
u(-\infty)=0\\
u(+\infty)=1.
\end{array}
\right. 
\end{equation}
Before going to the proof of Theorem
\ref{mono.th.sup}, we prove some \textit{a-priori} estimates on
$c^*_n$ . Namely we have the following,

\begin{proposition}{\ \  }\\
Let $c^*_n$ be the critical speed defined above, then there exists
a positive constant  $\kappa_1$ such that
$$-\kappa_1\le c^*_n.$$
\label{mono.prop.esti}
\end{proposition}
\dem{Proof.} Let  $f_\theta$ be a fixed function of ignition type
such that $f_\theta\le f$. Using Theorem \ref{mono.th.ig}, we have
$c^\theta_n\le c_n^*$. To obtain our desired bound, we just have
to bound from below  $c^\theta_n$. The later is obtained using
Lemma \ref{mono.lem.ig}. Indeed, for each $n\in\N$, we have
$$
\mu (c_n^{\theta})^2-\nu_n|c_n^\theta|\le 0,
$$
with $\nu_n:=\int_{\R}J_n(z)|z|\,dz$ and $\mu$ is independent of
$n$. Since $\nu_n\le \bar \nu:=\sup_{n\in\N}\{\nu_n\}<\infty$, we
end up with
$$
\mu (c_n^{\theta})^2-\bar \nu|c_n^\theta|\le 0.
$$
Hence,
$$
|c_n^{\theta}|\le \kappa_1.
$$
\fdem

Let us also recall some properties of  the
following  approximated problem:
\begin{equation}
\label{mono.eq.5.4} 
\left \{
\begin{array}{ll}
J_n\star u - u - cu^{\prime} +f\eta_{_{\theta}}(u)= 0 & \text{in} \  \ \R, \\
u(-\infty)=0,\\
u(+\infty)=1,
\end{array}
\right. 
\end{equation}
where $\theta>0$ and $\eta_{_{\theta}}$
is such that
\begin{itemize}
\item
$\eta_{_{\theta}}\in C_{0}^{\infty}(\R)$,
\item
$0\le \eta_{_{\theta}}\le 1$,
\item
$ \eta_{_{\theta}}(s)\equiv 0$ for $s\le \theta$ and  $
\eta_{_{\theta}}(s)\equiv  1$ for $s\ge 2\theta$.
\end{itemize}
For such kind of problem we have,

\begin{lemma} \label{lema 4.1}
Let $c^{\theta}$ and $c_n^{\theta}$ be the unique speed solution
of  \eqref{mono.eq.pre epsilon 0} with the nonlinearity $f
\eta_{_{\theta}}$ and respectively the kernel $J$ and $J_n$ and
let   $c_n^*$ be  the critical  speed solution of
\eqref{mono.eq.5.4} with  the nonlinearity $f$ and the kernel
$J_n$. Then the following holds:
\begin{itemize}
\item[a)] For fixed $\theta$, \quad
$\lim_{n\to \infty}c_n^{\theta}=c^{\theta}.$
\item[b)] For a fixed $n$, then \quad $\lim_{\theta\to 0}c_n^{\theta}=c_n^*.$
\end{itemize}
\label{mono.lem.limcn}
\end{lemma}
Part b) of this Lemma is contained in Proposition~\ref{mono.th.m}.
Part a) can be proved using similar arguments as in
Proposition~\ref{mono.th.m}.
\medskip


\subsection{\bf Construction of the solutions: Proof of Theorem \ref{mono.th.sup} }

We are now in position to prove Theorem \ref{mono.th.sup}. From
Proposition~\ref{mono.th.m}, we already know that there exists a
travelling front to \eqref{mono.eq.tw}-\eqref{mono.eq.twbc+} with
a critical speed $c^*$. To complete the proof, we have to
construct non decreasing solution for any speed $c\ge c^*$. We
emphasize that since $(w,\kappa)$ is not a super-solution of
\eqref{mono.eq.tw}-\eqref{mono.eq.twbc+} with the approximated
kernel $J_n$, there is   no uniform upper bound directly available
for the speed $c_n^*$ and the argumentation in the above section
cannot directly be applied.

From Proposition \ref{mono.prop.esti}, we have the following
dichotomy: either $\liminf (c_n^*)_{n\in\N}<+\infty$ or $\liminf
(c_n^*)_{n\in\N}=+\infty.$ We prove that in both situations there
exists a front solution for any speed $c\ge c^*$.

\bigskip

\noindent \underline{Case 1: $\liminf (c_n^*)_{n\in\N}<+\infty$}
\medskip

In this case, the same argument as in Proposition~\ref{mono.th.jc} in
Section~\ref{mono.ss.step3} works. Indeed, up to a subsequence
$c_n^* \to \tilde c$ and we must have $c^* \le \tilde c$. To prove
that $c^* = \tilde c$ we proceed as in
Section~\ref{mono.ss.step3}, using Lemma~\ref{lema 4.1} instead of
Lemma~\ref{limits}.

\fdem

Let now turn our attention to the other situation.
\medskip

\noindent \underline{Case 2: $\liminf (c_n^*)_{n\in\N}=+\infty$}
\medskip

In this case $\lim_{n\to\infty} c_n^*=+\infty$ we argue as
follows. Fix $c>c^* $, $c\not= 0$  where $c^*$ is defined by
Proposition~\ref{mono.th.m}. We will show that for such $c$ there
is a monotone solution to
\eqref{mono.eq.tw}-\eqref{mono.eq.twbc+}. When $c^* \le 0$ and
$c=0$ then a standard limiting procedure will show that a monotone
solution exists with this speed.

Again, by Theorem \ref{mono.th.ig} and
Proposition~\ref{mono.th.m}, we have $c^{\theta}<c^*$ for every
positive $\theta$. Therefore,
$$
\forall \theta>0,\; c^{\theta}< c^{*}< c.
$$
Fix $\theta > 0$. Since $c_n^{\theta}\to c^{\theta}$, one has on
one hand $c_n^{\theta}<c$ for $n\ge n_0$ for some integer $n_{0}$.
On the other hand,  $c^*_n\to +\infty$, thus there exists an
integer $n_1$ such that $c<c_n^*$ for all $n\ge n_1$. Therefore,
we may achieve for  $n\ge \sup\{n_0,n_1\}$,
$$
c_n^{\theta} <c<c_n^{*}.
$$

From this last  inequality, and according to Theorem
\ref{mono.th.ig} and Lemma \ref{mono.lem.limcn}, for each $n\ge
\sup\{n_0,n_1\}$ there exists a positive $\theta(n)\le \theta$
such that $c=c^{\theta(n)}_n$.

Let $u_n$ be the non decreasing solution of
\eqref{mono.eq.japprox} associated with $\theta(n)$. Since $
\theta(n)$ is bounded, we can extract a subsequence still denoted
$(\theta(n))_{n\in\N}$ which converges to some $\bar \theta$. We
claim that

\medskip\noindent{\bf Claim. $\bar \theta =0$}

\medskip

Assume for the moment that the claim is proved. Using the
translation invariance, we may   assume that for all $n$, $u_n(0)
=\frac{1}{2}$. Using now that $u_n$ is uniformly bounded and
Helly's theorem, up to a subsequence $u_n \to u$ pointwise, where
$u$ is a solution of \eqref{mono.eq.tw}-\eqref{mono.eq.twbc+} with
speed $c$.

In this way we get a non trivial solution of
\eqref{mono.eq.tw}-\eqref{mono.eq.twbc+} for any speed $c \ge
c^*$. \fdem

Let us now turn our attention to the proof of the above claim.

\dem{Proof of the Claim.} We argue by contradiction. If not, then
$\bar \theta >0$ and the speed $c{^{\bar \theta}}$ of the
corresponding non decreasing front solution of \eqref{mono.eq.pre
epsilon 0} satisfies
$$c^{{\bar\theta}}<c^*<c.$$
Let now consider, $u_n$ the solution associated with $\theta(n)$,
normalized by $u_n(0)=\theta(n)$. Using  uniform {\it a priori}
estimates,  Helly's theorem  we can extract a converging sequence
of function and get a solution $u$ with speed $c$ of the following
\begin{align*}
&J\star u - u - cu^{\prime} +f_{\bar\theta}(u)= 0\quad  \text{in}
\  \ \R .
\end{align*}
Using the arguments developed in \cite[Section 5.1]{Co4} to prove
Theorem~1.2 of that reference, one can show that $u$ satisfies the
boundary conditions
$$
u(+\infty)=1, \quad u(-\infty)=0.
$$
According Proposition~\ref{mono.th.m}, we get the contradiction
$$c=c^{{\bar \theta}}<c^*<c.$$
Hence $\bar \theta =0.$ \fdem


\section{Characterization of the minimal speed and asymptotic behavior}
\label{velocity}

Throughout this section we will assume the hypotheses of
Theorem~\ref{thm exp decay}, namely $f$ satisfies \eqref{f1},
$f'(0)>0$, $f\in C^{1,\gamma}$ near $0$ and \eqref{holder}, and
$J$ satisfies \eqref{j1}, $J\in C^1$ and is compactly supported.

Let us consider the following equation
\begin{equation}\label{tw2}
\begin{array}{l}
 J\star u -u- c u' + f(u) = 0 \quad \hbox{in
$\R$}, \\
\lim_{x\to-\infty} u(x) = 0.
\end{array}
\end{equation}

We need to establish some estimates on bounded solutions of
\eqref{tw2} that we constantly use along this section.

\begin{lemma} \label{mono.lem.estil1}
Let $u$ be a no-negative bounded solution of \eqref{tw2}, then the
following holds:
\begin{itemize}
\item[\textit{(i)}]
$ \int_{y}^{x}\int_{\R}J(s-t)[u(t)-u(s)]\,dtds=
\int_{0}^1\int_{\R}J(-z)z[u(x+z\eta)-u(y+z\eta)]\,dzd\eta $
\item[\textit{(ii)}]
$f(u)\in L^{1}(\R)$,
\item[\textit{(iii)}]
$u, J\star u \in L^{1}(\R^-) $
\item[\textit{(iv)}]
$v(x):=\int_{-\infty}^{x}u(s)\,ds$ satisfies $v(x)\le K(1+|x|)$
for some positive $K$ and $ v(x)\in L^1(\R^-)$.
\end{itemize}
\end{lemma}
\dem{Proof.} We start with the proof of \textit{(i)}. Let
$(u_n)_n$ be a sequence of  smooth $(C^1)$ function which converge
pointwise to $u$. Using the Fundamental Theorem of Calculus and
Fubini's Theorem, we have
\begin{align*}
\int_{y}^{x}\int_{\R}J(s-t)[u_n(t)-u_n(s)]\,dtds&=
\int_{y}^{x}\int_{0}^1\int_{\R}J(-z)zu_n'(s+z\eta)\,dzd\eta ds
\\
& =\int_{0}^1\int_{\R}J(-z)z[u_n(x+z\eta)-u_n(y+z\eta)]\,dzd\eta
\end{align*}
Since $|J(-z)zu_n(y+\eta z)|\le K|J(-z)z| \in
L^{1}(\R\times[0,1])$ and $u_n$ converges pointwise to $u$,
passing to the limit in the above equation yields
$$
\int_{y}^{x}\int_{\R}J(s-t)[u(t)-u(s)]\,dtds=
\int_{0}^1\int_{\R}J(-z)z[u(x+z\eta)-u(y+z\eta)]\,dzd\eta.
$$

To obtain  \textit{(ii)}, we argue as follow.
Integrating \eqref{tw2} from $y$ to $x$, it follows that
\begin{equation}
c(u(x)-u(y))-\int_{y}^{x}\int_{\R}J(s-t)[u(t)-u(s)]\,dtds
=\int_{y}^{x}f(u(s))\,ds. \label{mono.eq.asymp}
\end{equation}
Using \textit{(i)}, we end up with
\begin{equation}\label{boundintu}
c(u(x)-u(y))-\int_{0}^1\int_{\R}J(-z)z[u(x+z\eta)-u(y+z\eta)]\,dzd\eta
=\int_{y}^{x}f(u(s))\,ds.
\end{equation}
Again, since $|J(-z)zu(y+\eta z)|\le K|J(-z)z| \in
L^{1}(\R\times[0,1])$, we can pass to the limit  $y\to -\infty$ in
the above equation using Lebesgue dominated convergence Theorem.
Therefore, we end up with
$$
cu(x)-\int_{0}^1\int_{\R}J(-z)zu(x+z\eta)\,dzd\eta =\int_{-\infty}^{x}f(u(s))\,ds.
$$
Thus,
$$
\int_{-\infty}^{x}f(u(s))\,ds\le K( |c|+\int_{\R}J(z)|z|\,dz)
$$
which proves \textit{(ii)}. From equation \eqref{mono.eq.asymp},
we have
$$
c(u(x)-u(y))- \int_{y}^{x}f(u(s))\,ds + \int_{y}^{x}u(s)\,ds
=\int_{y}^{x} J\star u(s)\,ds.
$$
Thus $J\star u \in L^{1}(\R^-)$ will immediately follows from
$u\in L^{1}(\R^-)$ and \textit{(ii)}. Observe now that  since
$f'(0)>0$, and $u(-\infty)=0$, for $x<< -1$, we have  $f(u)>
\alpha u$ for some positive constant $\alpha $. Therefore,
$$
\alpha  \int_{-\infty}^{x}u(s)\,ds\le\int_{-\infty}^{x}f(u(s))\,ds
$$
and \textit{(iii)} is proved.

To obtain \textit{(iv)} we argue as follow. From \textit{(i-iii)},
$v$ is a well defined non decreasing function such that
$v(-\infty)=0$. Moreover, $v$ is smooth provide $u$ is continuous.
By definition of $v$, we easily see that $v(x)\le C(|x|+1)$ for
all $x\in \R$. Indeed, we have $$v(x)\le
\int_{\infty}^{0}u(s)\,ds+\int_{0}^{|x|}u(s)\,ds\le K(1+|x|),$$
where $K = \sup\{\int_{-\infty}^0u(s)\,ds; \nli{u}{\R}\}$.

Now, integrating \eqref{tw2} on $(-\infty,x)$, we easily see that
\begin{equation}\label{intform}
cv'(x)=J\star v(x)-v(x)+\int_{-\infty}^x f(u(s))ds.
\end{equation}
Since $f'(0)>0$ we can choose $R<<-1$ so that for $s\le R$,
$f(u(s))\ge \alpha  u(s)$ for some $\alpha >0$. Fixing  now $x<R$
and integrating \equ{intform} between $y$ and $x$, we obtain
\begin{equation}\label{aux2}
c(v(x)-v(y))\ge \int_y^x(J\star v(s)-v(s))ds+\alpha \int_y^x
v(s)ds.
\end{equation}
Proceeding as  above, we  get that $v\in L^1(-\infty,R)$.

\fdem

Following the idea of Carr and Chmaj \cite{CC}, we now derive some
asymptotic behavior of  the non negative bounded  solution $u$ of
\eqref{tw2} More precisely, we show the following

\begin{lemma}\label{lemmaexpdecayv}
Let $u$ be a non negative bounded continuous solution  of
\equ{tw2}. Then there exists two positive constants $M,\beta$,
such that $v(x)=\int_{-\infty}^x u(s) ds$ satisfies:
\begin{equation}\label{cotaexp u}
 v(x)\le Me^{\beta x},
 \end{equation}
\end{lemma}
\dem{Proof.} The proof  uses ideas from \cite{CD2}. Let first show
that some positive constants $C, R$,  we have
\begin{equation}\label{expdecay v}
\int_{-\infty}^{-R} v(x)e^{-\beta x}dx <C,
\end{equation}
for some $\beta>0$ small.

Consider $R>0$ and $\beta>0$ constants to be chosen later. Let
$\zeta\in C^\infty(\R)$ be a non-negative non-decreasing function
such that $\zeta\equiv0$ in $(-\infty,-2]$ and $\zeta\equiv1$ in
$[-1,\infty)$. For $N\in\N$, let $\zeta_N=\zeta(x/N)$. Multiplying
\eqref{intform} by $e^{-\beta x}\zeta_N$ and integrating over
$\R$, we get

\begin{equation}\label{mono.eq.zetan}
 \int_{\R} (J\star v-v)(e^{-\beta x}\zeta_N)\,dx - \int_{\R} cv'(e^{-\beta
x}\zeta_N)\,dx +\int_{\R}\int_{-\infty}^{x} f(u(s))\,ds(e^{-\beta x}\zeta_N) \,dx=0
\end{equation}

Note that by the monotonicity of $\zeta_N$ we have
\begin{align*}
\int_\R J\star v(x)\zeta_N(x)e^{-\beta x}\,dx&=\int_\R\int_\R J(x-y)e^{-\beta x}\zeta_N(x)v(y)dzdy \\
&=\int_\R\int_\R J(z)e^{-\beta(z+y)}\zeta_N(z+y)v(y)dzdy \\
&\ge \int_\R v(y)e^{-\beta y}\left(\int_{-R}^\infty J(z)e^{-\beta z}\zeta_N(
y-R)dz\right)\,dy.
\end{align*}
Therefore, we have
\begin{equation}
\int_{\R} (J\star v-v)(e^{-\beta x}\zeta_N)\,dx \ge \int v(x)e^{-\beta x}\left(\int_{-R}^\infty J(z)e^{-\beta
z}\;dz\;\zeta_N(x-R) - \zeta_N(x)\right)\;dx, \label{aux3}
\end{equation}

Let us now choose our adequate  $R>0$.  First pick
$0<\alpha<f'(0)$ and $R>0$ so large that
\begin{equation}
f(u)(x)\ge \alpha u(x) \qquad\text{for $x\le -R$.} \label{mono.eq.ineg1}
\end{equation}

Next, one can increase $R$ further if necessary so that
$\int_{-R}^\infty J(y)\;dy > (1-\alpha/2).$ By continuity we obtain
for some $\beta_0>0$ and all $0<\beta<\beta_0$,
\begin{equation} 
\label{mono.eq.ineg2}
\int_{-R}^\infty J(y)e^{-\beta y}\;dy \ge (1-\alpha/2)e^{\beta
R}. 
\end{equation}
Collecting \eqref{aux3} and \eqref{mono.eq.ineg2}, we then obtain

\begin{align}
\int_{\R} (J\star v-v)(e^{-\beta x}\zeta_N) &\ge \int_{\R} v(x)e^{-\beta
x}\left((1-\alpha/2)e^{\beta R}\zeta_N(x-R) -
\zeta_N(x)\right)\;dx\nonumber\\
&\ge (1-\alpha/2)\int_{\R} v(x+R)e^{-\beta x}\zeta_N(x)\;dx - \int_{\R} v(x)e^{-\beta x}\zeta_N(x)\;dx\nonumber\\
&\ge -\alpha/2\int_{\R} v(x)e^{-\beta x}\zeta_N(x)\;dx, \label{mono.eq.ineg3}
\end{align}
where we used the monotone behavior of $v$ in the last
inequality.

We now estimate the second term in \eqref{mono.eq.zetan}:

\begin{align}
\int_{\R} v'\zeta_N e^{-\beta x}\; dx &= \beta\int_{\R} v\zeta_N
e^{-\beta x}\; - \int_{\R} v \zeta_n'e^{-\beta x}\;dx\nonumber\\
&\le \beta\int_{\R} v\zeta_N e^{-\beta x}\;\label{mono.eq.ineg4}.
\end{align}
Finally using \eqref{mono.eq.ineg1}, the last term in \eqref{mono.eq.zetan}
satisfies

\begin{align}
\int_{\R}\left(\int_{-\infty}^{x} f(u(s))\,ds\right)\zeta_N e^{-\beta x}\;dx&=  \int^{-R}_{-\infty}
\left( \int_{-\infty}^{x} f(u(s))\,ds\right)\zeta_N e^{-\beta x}\;dx -C \nonumber\\
&\ge \alpha \int^{-R}_{-\infty}
v\zeta_N e^{-\beta x}\;dx -C. \label{mono.eq.ineg5}
\end{align}

By \eqref{mono.eq.zetan}, \eqref{mono.eq.ineg3}, \eqref{mono.eq.ineg4} and
\eqref{mono.eq.ineg5}, we then obtain

\begin{align*}
&|c|\beta \int_{\R} u\zeta_N e^{-\beta x}\;dx\ge  \alpha \int^{-R}_{-\infty} u\zeta_N e^{-\beta x}\;dx - C - \alpha/2   \int_{\R} v\zeta_N e^{-\beta x}\;dx \\
&(\alpha/2 -|c|\beta) \int^{-R}_{-\infty} u\zeta_N e^{-\beta x}\;dx \le \tilde C.
\end{align*}
Choosing $\beta <\alpha/(2|c|)$ and letting $N\to\infty$ proves \eqref{expdecay v}.

Using the monotonicity of $v$ we can conclude that
\begin{equation}\label{cotaexpv}
v(x)\le Ce^{\beta x},
\end{equation}
for some constant $C$. Indeed, if \equ{cotaexpv} does not hold,
then for a sequence $x_n\to-\infty$ we have $v(x_n)\ge n e^{\beta
x_n}$. Extracting a subsequence if necessary, we can assume that
$x_{n+1}<x_n -1$, thus since $v$ is increasing we have

\begin{align*}
\int_{-\infty}^{x_0} v(x) e^{-\beta x}dx&\ge \sum_{n\ge 1}\int_{x_n}^{x_{n-1}} n e^{\beta x_n} e^{-\beta x} dx\\
&\ge \sum_{n\ge 1}n\frac{1-e^{-\beta(x_n-x_{n-1})}}{\beta} \\
& \ge \sum_{n\ge 1}n\frac{1-e^{-\beta}}{\beta} =\infty
\end{align*}
which is a contradiction. \qed

In the next result we establish that the bounded solution $u$ of
\equ{tw2} also decays exponentially as $x\to -\infty$.

\begin{lemma}\label{expdecayu}
Suppose that $u$ is bounded  solution of \equ{tw2}. If for some $M,
\beta>0$ we have that $v(x)\le Me^{\beta x}$ for all $x$ then
 there exists $M_1,\alpha>0$ such that
\begin{equation}\label{cotaualpha}
u(x)\le M_1 e^{\alpha x} \ \hbox{for all}  \ x\in \R.
\end{equation}
\end{lemma}
\dem{Proof.} When $c\ne 0$ then  by \eqref{intform} we have the
following estimates
\begin{align*}
|c|u(x)&= \left| J\star v - v +\int_{\infty}^xf(u(s))\,ds\right|\\
          &\le  J\star v + v +\int_{\infty}^x\frac{f(u(s))}{u(s)}u(s)\,ds\\
           &\le  J\star v + (K+1)v
\end{align*}
where $K$ is the Lipschitz constant of $f$. Now since  $$J\star
v(x)\le C\int_{\R}J(x-y)e^{\beta y}\le C' e^{\beta x},$$ we easily
see that \equ{cotaexp u} holds.

When  $c=0$ the estimate does not directly comes from
\eqref{intform} and we have to distinguish several cases.

Let first observe that   for $x<0$ since $u$ is bounded by some
constant $C$,  $J\star u$  satisfies the following
\begin{align*}
J\star u(x)&=\int_{\infty}^{\frac{\alpha}{\beta}x}J(x-y)u(y)\,dy +
\int_{\frac{\alpha}{\beta}x}^{+\infty}J(x-y)u(y)\,dy
\\
&\le ||J||_\infty\int_{-\infty}^{\frac{\alpha}{\beta}x}u(y)
dy+C\int_{x\left(\frac{\alpha}{\beta}-1\right)}^\infty J(-z)\,
dz
\\
&\le ||J||_{\infty}v(\frac{\alpha}{\beta}x)+Ce^{(\beta-\alpha)
x}\int_{x\left(\frac{\alpha}{\beta}-1\right)}^\infty J(-z)e^{\beta
z}\,dz
\\
\end{align*}
Choosing $\alpha=\frac{\beta}{2}$ in the above equation,  we end
up with 
\begin{equation} 
\label{cotaJstaru} 
J\star u(x)\le C e^{\frac{\beta}{2}x},
\end{equation}
for some constant $C$. Observe also that since $f$ is smooth and $f(0)=0$, we have for small $\eps>0$   and $s>0$ small,
$$
|\frac{f(s)}{s}-f'(0)|\le \eps.
$$
Therefore from  \eqref{tw2}, for $\eps>$ small there exists
$K(\eps)>0$ such that for $x<-K(\eps)$ we have 
\begin{equation}
\label{mono.eq.expu} 
u(1-f'(0)+\eps)  \ge J\star u= u(1-\frac{f(u)}{u})\ge
u(1-f'(0)-\eps).  
\end{equation}
Observe now that if $f'(0)>1$, we get a contradiction. Indeed, choose $\eps$ so that $(1-f'(0)+\eps)<0$, then we have the following contradiction when $x<-K(\eps)$
$$
0>u(1-f'(0)+\eps)  \ge J\star u\ge 0.
$$
Thus, when $f'(0)>1$, there is no positive solution of \eqref{tw2}
with zero speed.

Let us now look at the other cases. Assume now  that $f'(0)<1$ and
choose $\eps$ small so that $(1-f'(0)-\eps)>0$ then from
\eqref{mono.eq.expu} for $x<-K(\eps)$ there exists a positive
constant $C$ so that
$$ u\le C J\star u\le C e^{\frac{\beta}{2} x}.$$
Finally, when $f'(0)=1$ recall that  $f $ satisfies
\eqref{holder}. Thus, for $x<<-1$
$$
J\star u(x)= u-f(u) \ge A u^{m},
$$
where $A>0$, $m\ge 1$. Using \eqref{cotaJstaru}, yields
$$
u\le \frac{C}{A} e^{\frac{\beta}{2 m} x}.
$$
\fdem

\begin{remark}\label{remark exp behavior}
From the above proof, we easily  conclude that for any $0<\alpha<
\bar \alpha $, where $\bar \alpha$ depends only on $\beta$ and
$\gamma$, there exists $M_1>0$ such that \equ{cotaualpha} holds.
\end{remark}
\medskip

As in \cite{CC}, for $u$  a solution of \equ{tw2}  we define the
function $U(\lambda)=\int_{\R} e^{-\lambda x} u(x)dx$ which by
Lemma \ref{expdecayu} is defined and analytic in the strip $0< Re
\, \lambda<\alpha$.  Note that
$$
\int_{\R} J\star u(x) e^{-\lambda x}=\int_\R u(y)e^{-\lambda
y}dy\int_\R J(-z)e^{\lambda z} dz
$$
and using integration by parts
$$c\int_{\R}u'e^{-\lambda x} dx=\lambda c\int_{\R}u(y)e^{-\lambda y}dy.$$
Using the above identities, if we multiply \equ{tw2} by
$e^{-\lambda x}$ and integrate in $\R$ we obtain
\begin{equation}
U(\lambda)(-c\lambda +m(\lambda))=\int_{\R} e^{-\lambda
x}(f'(0)u(x)-f(u(x)))dx, \label{eqUlambda}
\end{equation}
where the function $m(\lambda)=\int_{\R} J(-x)e^{-\lambda
x}dx+f'(0)-1$  is analytic in $\C$.

Let  $c^1$ be the following quantity
$$
c^1:= \min_{\lambda>0} \frac{1}{\lambda} \left( \int_\R J(-x)
e^{\lambda x} \, d x + f'(0) -1 \right).
$$
\begin{proposition}\label{propnonex}
If $c<c^1$ then \equ{tw2} does not have any solution.
\end{proposition}

\dem{Proof.} Since  $u>0$ we deduce, from a property of Laplace
transform (Theorem 5b, p. 58 \cite{widder}) and Lemma
\ref{expdecayu},  that the function $U(\lambda)$ is analytic in
$0<Re \,\lambda<B$, where $B\ge \alpha$,  and $U(\lambda)$ has a
singularity at $\lambda=B$. Observe that if $c<c^1$ then for some
$\delta>0$
\begin{equation}
\label{cotadelta} -c\lambda +m(\lambda)>\delta \lambda, \ \ \hbox{
for all} \ \lambda>0.
\end{equation}

Observe that since $f \in C^{1,\gamma} $ near 0  and using Lemma
\ref{expdecayu} we have that for some constant $C>0$

\begin{align*}
\int_{\R} e^{-\lambda x}|f'(0)u(x)-f(u(x))|dx&=\int_{\infty}^{-K}
e^{-\lambda x}|f'(0)u(x)-f(u(x))|dx
\\
& \qquad +  \int_{-K}^{+\infty} e^{-\lambda
x}|f'(0)u(x)-f(u(x))|dx
\\
&\le \int_{\infty}^{-K} e^{-\lambda x}|A
u^{1+\gamma}+o(u^{1+\gamma})|dx  +  C\int_{-K}^{+\infty}
e^{-\lambda x}u (x)dx
\\
&\le C\int_{\R}  e^{-\lambda x} u^{1+\gamma}(x) dx\\
&\le C\int_{\R}  e^{(-\lambda +\gamma \alpha)x} u(x) dx.
\end{align*}

From the above computation, it follows that $\int_{\R} e^{-\lambda
x}|f'(0)u(x)-f(u(x))|dx$ is analytic in the region $0<Re \,\lambda
<B + \gamma \alpha$. Since $\gamma>0$, using the equation
\equ{eqUlambda}, we get $U(\lambda)$ defined and analytic for
$0<Re \lambda<B+\gamma \alpha$. Bootstrapping this argumentation
we can extend analytically $U(\lambda)$ to $Re \, \lambda>0$. Then
for all $\lambda
>0$
$$
\int_{\R} e^{-\lambda x}|f'(0)u(x)-f(u(x))|dx\le (f'(0)+k)
\int_{\R} e^{-\lambda x}u(x)=CU(\lambda).
$$
Therefore for all $\lambda>0$, using \eqref{eqUlambda}, it follows
that $-c\lambda+m(\lambda)\le C$ contradicting \equ{cotadelta}.
\fdem

\begin{remark}\label{remark weaker}
We should point out that the above proposition holds as well if
the kernel $J$ instead of being compactly supported, is only
assumed to satisfy:
$$\exists\; M,\lambda_0>0\quad \text{ such that }\quad
\int_{0}^{+\infty}J(-x)e^{\lambda_0 x}\le M.$$
\end{remark}
\medskip

Let us now  establish the exact asymptotic behavior, as
$x\to-\infty$, of a solution $u$ of \equ{tw2}. We proceed as
follows. First, we  obtain the exact behavior of
$v=\int_{-\infty}^x u(s) ds$, proceeding as in \cite{CC} and then
we conclude the behavior of $u$.

For $c\ge c^1$ we  denote $\lambda(c)$ the unique minimal
$\lambda>0$ such that $-c\lambda +m(\lambda)=0$. It can be easily
verified that $\lambda(c)$ is a simple root of $-c\lambda
+m(\lambda)$ if $c>c^1$, and it is a double root when $c=c^1$.

\dem{Proof of Theorem~\ref{thm exp decay}.} Since there is a
monotone solution $(u,c^*)$ of
\eqref{mono.eq.tw}-\eqref{mono.eq.twbc+} with critical speed, it
is a bounded solution of \eqref{tw2}. Thus by
Proposition~\ref{propnonex} $c^*\ge c^1$.

It remains to prove \eqref{behavior1} and \eqref{behavior2}. The
proof follows from a modified version of Ikehara's Theorem (see
\cite{widder}). We define $F(\lambda)=\int_{-\infty}^0 v(y)
e^{-\lambda y}$. Since $v$ is monotone, we can obtain the
appropriate asymptotic behavior of $v$ if $F$ has the
representation
\begin{equation}\label{formula1}
F(\lambda)=\frac{H(\lambda)}{(\lambda-\alpha)^{k+1}},
\end{equation}
with $H$ analytic in the strip $0<  Re \lambda \le \alpha$, and
$k=0$ when $c>c^*$, $k=1$ when $c=c^*$.

Using \equ{intform}, we have that
$$
\int_{-\infty}^0 v(x) e^{-\lambda x} dx=
\frac{\int_{-\infty}^{\infty }\int_{-\infty}^x f(u(s))-f'(0)u(s)
\,  ds e^{-\lambda x} \, d x }{c\lambda-m(\lambda)}
-\int_0^{\infty} v(x) e^{-\lambda x},
$$
thus, using that either $c\ne 0$ or $f'(0)<1$ holds, we have that
by Lemma \ref{expdecayu}, \equ{formula1} holds replacing $u$ by
$v$ with $\alpha=\lambda(c)$ described above, since it can be
checked that $-c\lambda+m(\lambda)$ has only two real roots which
are simple when $c>c^1$ and double when $c=c^1$

It remains to conclude that \equ{formula1} holds for $u$. First
suppose  that $c=c^1$ and denote $\lambda=\lambda(c^1)$. If $c\ne
0$ then using \equ{intform} we have that
$$
cu=J\star v(x)-(1-f'(0))v(x)+\int_{-\infty}^x f(u(s))-f'(0)u(s) ds,
$$
By Remark \ref{remark exp behavior} and since $f$ is $C^{1,\gamma}$ near 0 we have that
\begin{equation}\label{limite auxiliar1}
\frac{\int_{-\infty}^x f(u(s))-f'(0)u(s) ds}{ |x |e^{-\lambda
(c^1)x}}\to 0,
\end{equation}
as $|x|\to -\infty$. Therefore, we just have to prove that
\begin{equation}
\label{limiting exponential} \lim_{x\to-\infty} \frac{J\star
v(x)-(1-f'(0))v(x)}{ |x| e^{\lambda (c^1) x}}=L\ne 0.
\end{equation}
Observe that since $v$ satisfies \equ{behavior1} we have that for
${\eta}=\lim_{x\to-\infty} \frac{v(x)}{ |x| e^{\lambda (c^1) x}}$
and $supp \, J \subset [-k,k]$ we have
$$
\displaystyle{\frac{J\star v}{|x| e^{\lambda(c^1)
x}}}=\frac{1}{|x|}\int_{-k}^k
J(-z)(\eta+O(1/x))e^{\lambda(c^1)z}(|x|+z) dz,
$$
therefore
$$
\frac{J\star v(x)-(1-f'(0))v(x)}{ |x| e^{\lambda (c^1) x}}\to \eta
m(\lambda(c^1))=\eta c^1\lambda(c^1)\ne 0,
$$
which gives the desired result.

When $c^1=0$, we proceed in a slightly different way. Observe that
in this case $f'(0)<1$
\begin{equation}
(1-f'(0)) u= J\star u +f(u)-f'(0)u, \label{despeje}
\end{equation}
and by Remark \ref{remark exp behavior} and since $f \in C^{1,\gamma}$ near 0 we have
that \equ{limite auxiliar1} holds. Also, by \eqref{j2} we have
that $J\star u=J'\star v$ and
$$
\begin{array}{rl}
\displaystyle{\frac{J'\star v}{|x| e^{\lambda(c^1)
x}}}&=\frac{1}{|x|}\int_{-k}^k
J'(-z)(\eta+O(1/x))e^{\lambda(c^1)z}(|x|+z) dz
\\
\\
&=\eta \int _\R J(-z) e^{\lambda(c^1)z} dz + O(1/x),
\end{array}
$$
with $\eta >0$ as above. Hence, we obtain the desired result.

Finally, the case $c>c^1$ is analogous. \fdem

\dem{Proof of Corollary~\ref{corollary velocity kpp}.} Observe now
that in the case of a KPP nonlinearity $f$, the function
$w:=e^{\lambda x}$ is a super-solution of
\eqref{mono.eq.tw}-\eqref{mono.eq.twbc+}, provided that
$\lambda>0$ is chosen such that $-c\lambda +m(\lambda)=0$. The
existence of such $\lambda>0$ is guaranteed since $c\ge c^1$. The
existence of a monotone travelling wave for any $c\ge c^1$ is then
provided by Theorem \ref{mono.th.sup}. Therefore $c^*\le c^1$ and
we conclude $c^* = c^1 $. \fdem


\section{Uniqueness of the profile}
\label{section uniqueness}

In this section we deal with the uniqueness up to translation of
solution of \eqref{mono.eq.tw}-- \eqref{mono.eq.twbc+}. Our proof
follows ideas of  \cite{Co2} and is mainly  based on the sliding
methods introduced by Berestycki and Nirenberg \cite{BN1,BN2} (see also \cite{Co2}).

In the sequel, given a function $u:\R\to\R$ and $\tau \in \R$ we
define its translation by $\tau$ as
\begin{align} \label{notation translation}
u_\tau(x) = u(x+\tau)
\end{align}
and sometimes we shall write $u^\tau(x)= u(\tau+x)$.

Let $L$ denote the operator
$$
L u = J \star u - u - c u'.
$$

\begin{proposition} \label{Nonlinear Comparison Principle}
(Nonlinear Comparison Principle)\\
Let $J$ satisfy  \eqref{j1}, \eqref{MP} and let $f$ be a
monostable nonlinearity so that $f'(1)<0$. Let  $u$ and $v$ be two
continuous functions in $\R$ such that
\begin{align}
  &Lu + f(u)\le 0 \ \ \text{ on }\ \ \R \label{mocovequ}\\
  &Lv + f(v)\ge 0 \ \ \text{ on }\ \ \R \label{mocoveqv}\\
  & \lim_{x\to -\infty} u(x)\ge 0, \ \ \lim_{x\to -\infty} v(x)\le 0 \ \ \label{mocovbc-}\\
  & \lim_{x\to +\infty} u(x)\ge 1,  \ \ \lim_{x\to +\infty} v(x)\le 1.\ \  \label{mocovbc+}
\end{align}
Assume further that either  $u$ or $v$ is monotone and that $u \ge
v$ in some interval  $(-\infty,K)$. Then there exists $\tau \in
\R$ such that $u_\tau \ge v$ in $\R$. Moreover, either $u_\tau> v$
in $\R$ or $u_\tau\equiv v$.\label{mothmcomp}
\end{proposition}

\begin{remark}
Observe that by the Maximum Principle  and since $f(s)\ge 0\ \
\forall s\le 0$, the supersolution $u$ is necessarily positive.
Similarly, since $f(s)\le 0\ \ \forall s\ge 1$, the Maximum
Principle implies that $v<1$.  \label{mocovrem}
\end{remark}

\dem{Proof of Proposition~\ref{mothmcomp}.} Note that if
$\inf_{\R} u \ge \sup_{\R} v$, the theorem trivially holds. In the
sequel, we assume that  $\inf_{\R} u < \sup_{\R} v$.

Let $\eps>0$ be such that
\begin{equation} \label{mo1}
f^{\prime}(p)\le 0 \quad \text{ for } \ \  1-\eps<p<1.
\end{equation}
Now fix $0<\delta\le \frac{\eps}{2}$  and choose $ M>0$
sufficiently large so that
\begin{align}
   &\ \ 1-u(x) <\frac {\delta}{2} \quad  \forall x >M \label{mo2}
   \\
   &\ \ v(x)< \frac {\delta}{2} \quad \forall x <-M \label{mo3}
   \\
   \text{ and } &\ \ v(x) \le u(x)  \quad \forall x <-M. \label{mo4}
\end{align}

\medskip
\noindent{\bf Step 1.} There exists a constant $D$ such that for
every $b\ge D$
\begin{align} \label{a1}
u(x+b)>v(x) \quad \forall x \in [-M-1-b,M+1].
\end{align}
Indeed, since $u>0$ in $\R$ and $\lim_{x\to+\infty} u(x) \ge 1$ we
have
$$
c_0 := \inf_{[-M-1,\infty)} u >0
$$
Since $\lim_{x\to-\infty} v(x) \le 0$ there is $L>0$ large such
that
$$
v(x) < c_0 \quad \forall x\le -L.
$$
Then for all $b>0$
$$
u(x+b) > v(x) \quad \forall x \in[ M-1-b,-L].
$$
Now, since $\sup_{[-L,M+1]} v < 1$ and $\lim_{x\to+\infty} u(x)
\ge 1$ we deduce \eqref{a1}.

\bigskip \noindent{\bf Step 2.} There exists $b\ge D$ such that
\begin{align} \label{a2}
u(x+b)+\frac{\delta}{2}>v(x)\ \ \forall x \in \R .
\end{align}
If not then we have,
\begin{align} \label{moneg}
\forall b\geq D\  \ \text{ there exists}\ \ x(b)\ \ \text{such
that}\ \ u(x(b)+b)+\frac{\delta}{2}\leq v(x(b)).
\end{align}
Since $u$ is nonnegative and  $v$   satisfies \eqref{mocovbc-}
there exists a positive constant A such that
\begin{align} \label{mo6}
u(x+b)+\frac{\delta}{2}> v(x) \ \ \text{ for all }  b>0 \text{ and
}\ \ x \leq -A.
\end{align}
Take now a sequence $(b_n)_{n\in\N}$ which tends to $+\infty$. Let
$x(b_n)$ be the point defined by \eqref{moneg}. Thus we have for
that sequence
\begin{align} \label{mo7}
u(x(b_n)+b_n)+\frac{\delta}{2}\leq v(x(b_n)).
\end{align}
According to \eqref{mo6} we have $x(b_n)\geq -A.$ Therefore the
sequence $x(b_n)+b_n$ converges to $+\infty$. Pass to the limit in
\eqref{mo7} to get
\begin{equation*}
1+\frac{\delta}{2} \le \lim_{n\to +\infty}
u(x(b_n)+b_n)+\frac{\delta}{2}\le  \limsup_{n\to +\infty}
v(x(b_n)) \le 1,
\end{equation*}
which is a contradiction. This proves our claim \eqref{a2}.

\bigskip \noindent{\bf Step. 3}
We observe that as a consequence of \eqref{a1} and \eqref{a2}, and
using that either $u$ or $v$ is monotone we in fact have
\begin{align}
u(x+b) & \ge v(x) \quad \forall x \le M+1 \label{a4}
\\
u(x+b) + \frac{\delta}{2}& > v(x) \quad \forall x \ge M+1 .
\nonumber
\end{align}
Indeed, it only remains to verify that $u(x+b)>v(x)$ for $x \le
M-1-b$. If $u$ is monotone from \eqref{mo4} we have
$u(x+b)>u(x)>v(x)$ for $x<-M$. If $v$ is monotone $u(x) >
v(x)>v(x-b)$ for $x<-M$.

\bigskip \noindent{\bf Step 4.} Now we claim that
\begin{align} \label{a3}
u(x+b) \ge v(x) \quad \forall x \in \R.
\end{align}
To prove this, consider
\begin{align} \label{moa2}
a^*=\inf\{ a>0\mid  u(x+b)+a \ge v(x)\;\forall x \in \R\}
\end{align}
which is well defined by \eqref{a2}.

If $a^*= 0$ then \eqref{a3} follows. Suppose $a^*>0$. Then, since
$$
\lim_{x\to \pm \infty}u(x+b)+a^* -v(x)\ge a^*>0,
$$
there exists $x_0 \in \R$ such that   $u(x_0+b)+a^*=v(x_0)$.

Let $w(x):=  u(x+b)+a^*-v(x)$ and note that
\begin{align} \label{momin}
0=w(x_{0})=\min_{\R} w(x).
\end{align}
Observe  that $w$ also satisfies  the following equations:
\begin{align}
  &Lw \le  f(v(x))-f(u(x+b)) \label{moeqp}\\
  &w(+ \infty )\ge a^*  \label{mocl1} \\
  &w(- \infty )\ge a^*.  \label{mocl2}
\end{align}
Since $w\ge 0$, $ w \not\equiv 0$ using the strong maximum
principle at $x_0$ we have
\begin{align} \label{mowmin}
Lw(x_0)> 0.
\end{align}
By \eqref{a4} we necessarily have $x_0>M+1$.

At $x_0$ we have
\begin{align} \label{moq+}
f(u(x_0+b)+a^*)-f(u(x_0+b))\le 0,
\end{align}
since $f$ is non-increasing for $s\ge 1- \eps$, $a^*>0$ and
$1-\eps< 1- \frac{\delta}{2}\le u$ for $x>M$. Combining
\eqref{moeqp},\eqref{mowmin} and \eqref{moq+} yields  the
contradiction
$$0<Lw(x_0)\le f(u(x_0+b)+a^*)-f(u(x_0+b)) \le 0.$$

\medskip \noindent{\bf Step 5.}
Finally it remains to prove that either $u_\tau> v$ or $u_\tau
\equiv v$. Let $w:=u_\tau - v $, then either $w>0$  or $w(x_0)=0$
at some point $x_0 \in \R$. In the latter case we have $w(x)\ge
w(x_0)=0$ and
\begin{equation}
0\le Lw(x_0) \le  f(v(x_0))-f(u(x_0+\tau))=f(v(x_0))-f(v(x_0))= 0.
\end{equation}
Then using the maximum principle, we obtain $w\equiv 0$, which
means $u_\tau \equiv v$. \fdem

\bigskip

\begin{proposition} \label{mono.lem.order}
Let $J$ satisfy  \eqref{j1}, \eqref{MP} and let $f$ be a
monostable nonlinearity so that $f'(1)<0$. Let $u_1$ and $u_2$ be
respectively
 super  and sub-solutions of
\eqref{mono.eq.tw}-\eqref{mono.eq.twbc+} which are continuous. If
$u_1\ge u_2$ in some interval $(-\infty,K)$ and either $u_1$ or
$u_2$ is monotone then $u_1\ge u_2$ everywhere. Moreover either
$u_1>u_2$ or $u_1\equiv u_2$.
\end{proposition}

\dem{Proof.} Assume first that $\inf_{\R} u_1 < \sup_{\R} u_2$.
Otherwise there is nothing to prove. Without losing generality we
can assume that $u_1$ is monotonic. Using Theorem \ref{mothmcomp},
$u_1^\tau\ge u_2$ for some $\tau\in \R$, so  the following
quantity is well defined
$$
\tau^*:=\inf\{\tau\in \R| u_1^\tau\ge u_2\}
$$
We claim that
\begin{align} \label{a5}
\tau^*\le 0
\end{align}

Observe that by showing that  $\tau^*\le 0$,  we end the proof. To
prove \eqref{a5} we argue by contradiction. Assume that
$\tau^*>0$, then since $u_i$ are a continuous functions, we will
have $ u_1^{\tau^*} \ge u_2$ in $\R$. Let $w:=u_1^{\tau^*}-u_2\ge
0$. Since $\tau^*>0$ and $u_1$ is monotone then $w>0$ in
$(-\infty,K)$. Now observe that $w>0$ in $\R$ or $w(x_0)=0$ for
some point $x_0$ in $\R$. In the latter case
$$
0\le (J\star w-w)(x_0)\le f(u_2(x_0))-f(u_1^{\tau^*}(x_0))=0.
$$
Thus, using the maximum principle, $w\equiv 0$, which contradicts
that $w>0$ in $(-\infty,K)$. Now since $u_1$ is monotonic and
$\tau^*>0$ for small $\eps>0$, we have $u_1^{\tau^*-\eps}>u_2$ in
$(-\infty,M)$. Arguing as in Step 4 of the proof of
Proposition~\ref{mothmcomp} we deduce $u^{\tau^*-\eps}>u_2$ in
$\R$ which contradicts the definition of $\tau^*$. \fdem

\begin{remark} \label{one continuous}
With minor modifications the proofs of the
Propositions~\ref{Nonlinear Comparison Principle} and
\ref{mono.lem.order} hold if only one of the functions $u_1$ or
$u_2$ is continuous. For the proof of this statement we need the
strong maximum principle for solutions in $L^\infty$, which can be
found in \cite{CDM}:

{\rm \begin{theorem} \label{strong max p a.e.} Assume $J$
satisfies \eqref{j1}, \eqref{MP} and let $c \in L^\infty(\R)$. If
$u \in L^\infty(\R)$ satisfies $u \le 0 $ a.e.\@ and $J \star u -
u  + c(x) u \ge 0$ a.e.\@ in $\R$, then $ess\ sup_K u <0 $ for all
compact $K\subset \R$  or $u =0 $ a.e.\@ in $\R$.
\end{theorem}}
\end{remark}

\dem{Proof of Theorem \ref{mono.th.uniq}.} The case of  $c\neq
c^1$ and $c=c^1$ being similar, we present only the case $c\neq
c^1$. Let $u_1$ and $u_2$ be two    solution of
\eqref{mono.eq.tw}-\eqref{mono.eq.twbc+} with the same speed
$c\neq 0$. Since $c\ne 0$ the functions  $u_i$ are uniformly
continuous. From Theorem \ref{mono.th.sup}, we can assume that
$u_1$ is a monotonic function. Since, $u_i$ solve the same
equation and $u_1$ is monotone, using the translation invariance
of the equation and \eqref{behavior2} we see that up to a
translation
\begin{align}
&u_1= e^{\lambda(c)x} + o(e^{\lambda(c)x})\label{mono.eq.adu1}
\quad \hbox{as $x\to-\infty$}
\\
&u_2= e^{\lambda(c)x} + o(e^{\lambda(c)x}) \quad \hbox{as
$x\to-\infty$}.\label{mono.eq.adu2}
\end{align}
Let us first recall the following notation,
$u^{\tau}(.):=u(.+\tau)$. Then, by monotonicity of  $u_1$ and
\eqref{mono.eq.adu1}- \eqref{mono.eq.adu2} for some positive
$\tau$ we have  $u_1^\tau\ge u_2$ in some interval $(-\infty,-K)$.
Using Proposition~\ref{mono.lem.order}, it follows that
$u_1^{\tau}\ge u_2$ for possibly a new $\tau$. Define now the
following quantity:
$$\tau^*:=\inf \{\tau>0| u_1^{\tau}\ge u_2\}$$
Observe that form the above argument $\tau^*$ is well defined. We
claim

\medskip\noindent{\bf Claim: $\tau^*=0$.}

Observe that proving the claim ends the proof of the uniqueness up
to translation of the solution. Indeed, assume for a moment that
the claim is proved then we end up with $u_1\ge u_2$.
Observe now that in the above argumentation  the role of $u_1$ and $u_2$ can be
interchanged, so we easily see that we have $u_1\le u_2\le u_1$ which
ends the proof of the uniqueness. \fdem

Let us now prove the Claim.

\dem{Proof of the Claim.} If not, then $\tau^*>0$. Let
$w:=u_1^{\tau^*}-u_2 \ge 0$. Then either there exists a point
$x_0$ where $w(x_0) = 0$ or $w>0$. In the first case, at $x_0$,
$w$ satisfies:
$$
0\le J\star w(x_0) -w(x_0) =f(u_2(x_0))-f(u_1^{\tau^*}(x_0))=0
$$
Using the strong maximum principle, it follows that $w\equiv 0$.
Thus $u_1^{\tau^*}\equiv u_2$, which contradicts
\eqref{mono.eq.adu1}--\eqref{mono.eq.adu2}. Therefore,
$u_1^{\tau^*}>u_2$. Using  \eqref{mono.eq.adu1}, since $\tau^*>0$
we have for $u_1^{\tau^*}$  the following behavior near $-\infty$.
$$
u_1^{\tau^*}:=e^{\tau^*}e^{\lambda(c)x} + o(e^{\lambda(c)x}).
$$
Therefore, for some $\eps>0$ small, we still have
$u_1^{\tau^*-\eps}\ge u_2$ in some neighborhood $(-\infty,-K)$ of
$-\infty$. Using Theorem \ref{mono.lem.order}, we end up with
$u_1^{\tau^*-\eps}\ge u_2$ everywhere, contradicting the
definition of $\tau^*$. Hence, $\tau^*=0$. \fdem

\bigskip

Regarding Theorem~\ref{theo-uniquenessCC} we need the following
result:

\begin{lemma} \label{lemma liminf positive}
Assume that $J$ and $f$ satisfy \eqref{j1}, \eqref{j2}, \eqref{MP}
and \eqref{f1}, \eqref{f2} respectively. Let $0 \le u \le 1$ be a
solution to \eqref{mono.eq.tw}.

a) Then
$$
\lim_{x\to -\infty} u(x)=0 \quad \hbox{or}\quad
\lim_{x\to-\infty}u(x) =1,
$$
and
$$
\lim_{x\to \infty} u(x)=0 \quad \hbox{or}\quad
\lim_{x\to\infty}u(x) =1.
$$

b) If $u(-\infty) = 1$ and $u(+\infty)=1$ then $u\equiv1$.
\end{lemma}

Note that in this lemma we do not assume that $u$ is continuous.

\dem{Proof.}

a) Let $0\le u\le1$ be a solution to \eqref{nonlocal ind of x}. We
first note that by \equ{boundintu} any bounded solution $u$ of
\eqref{mono.eq.tw} satisfies
\begin{align} \label{f(u) integrable}
\int_{-\infty}^\infty f(u) du<\infty.
\end{align}
Let $g(u) = u - f(u)$ and note that
\begin{align}\label{eq J g}
J\star u = g(u) \quad \hbox{in $\R$},
\end{align}
and that the hypotheses on $f$ imply $g'(u) \ge g'(0)$ and $ g(u)
\le u$ for $ u \in [0,1]$.

If $f'(0)<1$ then $g'(0)>0$ and then $g$ is strictly increasing.
This together with \eqref{eq J g} implies that $u$ is uniformly
continuous and using \eqref{f(u) integrable} we see that
$u(-\infty)=0$ or $u(-\infty)=1$ and the same at $+\infty$ which
is the desired conclusion. Therefore in the sequel we assume
$f'(0) \ge 1$, that is, $g'(0) \le 0$.

\medskip

Since both limits at $- \infty$ and $+ \infty$ are analogous we
concentrate on the case $x\to -\infty$.

We will establish the conclusion of part a) by proving
\begin{align}\label{alternative liminf lim}
\liminf_{x\to-\infty} J \star u (x) =0  \quad \Longrightarrow
\quad \lim_{x \to -\infty} u(x) =0,
\end{align}
and
\begin{align}\label{alternative liminf lim 2}
\liminf_{x\to-\infty} J \star u (x) >0 \quad \Longrightarrow \quad
\lim_{x \to -\infty} u(x) =1.
\end{align}

We start with \eqref{alternative liminf lim}. Suppose that
$f'(0)>1$. Then there is $\delta
>0$ such that $g(u)<0$ for $u \in (0,\delta)$ and from \eqref{eq J
g} we deduce that $u(x) \ge \delta $ for all $x$, so regarding
\eqref{alternative liminf lim} there is nothing to prove.

Suppose $f'(0)=1$. Then $g$ is non-decreasing and by
\eqref{holder} we have, for some $A>0$, $m\ge 1$, $\delta_1>0$
\begin{align}\label{holder g}
g(u) \ge A u^m \quad \forall 0 \le u \le \delta_1.
\end{align}
Assume that $\liminf_{x\to-\infty} J \star u (x) =0 $ and let us
show first that
\begin{align}\label{limit zero}
\lim_{x\to-\infty} J \star u (x) = 0.
\end{align}
Otherwise, set $ \overline l = \limsup_{x\to-\infty} J \star u (x)
>0$. Choose $l \in (0,\overline l)$ such that $g'(l)>0$ and then
pick a sequence $x_n \to -\infty$ such that $J \star u(x_n)=g(l)$
for all $n$. Then there is some $\sigma>0$ such that for
$x\in(x_n-\sigma,x_n+\sigma)$ we have $f(u(x)) \ge c >0$ for some
uniform $c$. This contradicts \eqref{f(u) integrable} and we
deduce \eqref{limit zero}. This combined with \eqref{holder g}
implies that $\lim_{x \to -\infty} u(x) =0$, and this establishes
\eqref{alternative liminf lim}.

We prove now \eqref{alternative liminf lim 2}. Let us assume
$$
\underline l := \liminf_{x\to-\infty} J \star u (x) > 0.
$$
Since $ J \star u = g(u) \le u$ it is enough to show that
\begin{align} \label{lim J=1}
\lim_{x\to-\infty} J \star u (x) =1.
\end{align}
Assume the contrary, that is
\begin{align} \label{liminf <1}
0< \underline l  < 1.
\end{align}
Observe that
$$
\liminf_{x\to-\infty} u(x) >0.
$$
This is direct if $f'(0)>1$ and follows from \eqref{eq J g},
\eqref{holder g} and $\underline l>0$ if $f'(0)=1$. Therefore
$\limsup_{x\to-\infty}u(x)=1$, otherwise \eqref{f(u) integrable}
can not hold. Hence
\begin{align}\label{limsup =1}
\limsup_{x\to-\infty} J \star u (x) = 1.
\end{align}
Chose now $\alpha \in ( \underline l , 1) $ a regular value of the
function $g$. By \eqref{liminf <1}, \eqref{limsup =1} and the
continuity of $J\star u$ there exists a sequence $x_n \to -
\infty$ such that $J\star u(x_n) = \alpha$. Note that the set $ \{
u \in [0,1] \, /\, g(u) = \alpha \}$ is discrete and hence finite
and does not contain $0$ nor $1$. Hence, for sufficiently small
$\eps>0$ we have $ \{ u \in [0,1] \, /\, \alpha-\eps< g(u) <
\alpha + \eps \} \subseteq [\eps,1-\eps]$. Since $J\star u$ is
uniformly continuous there is $\sigma>0$ such that for $x\in
(x_n-\sigma,x_n+\sigma)$ we have $\eps \le u(x_n) \le 1 - \eps$.
This contradicts the integrability condition \eqref{f(u)
integrable}, and we deduce the validity of \eqref{lim J=1}.

\bigskip

b) Assume that $\lim_{x\to \infty} u(x)=\lim_{x\to-\infty}u(x)=1$
and set $\gamma^*= \sup\{ 0<\gamma <1 \ / \ u
>\gamma\}$. For the sake of contradiction assume that $u$ is
nonconstant. Then $0<\gamma^* <1$. Since $f(\gamma^*)>0$ we have
that $v=u-\gamma^*\ge 0$ satisfies
\begin{equation}\label{auxlimits}
J\star v-v-cv'+\frac{f(u)-f(\gamma^*)}{u-\gamma^*} (u-\gamma^*)<0.
\end{equation}
If $c\not=0$ then $v$ reaches its global minimum at some
$x_0\in\R$ which satisfies $v(x_0)=0$. Thus, evaluating
\eqref{auxlimits} at $x_0$ we obtain a contradiction. If $c=0$ we
reach again a contradiction applying Theorem~\ref{strong max p
a.e.}.

\fdem

\bigskip

\dem{Proof of Theorem~\ref{theo-uniquenessCC}.} Assume $0\le u \le
1$ is a solution of \eqref{nonlocal ind of x} such that
$u\not\equiv 0$ and $u\not\equiv 1$. By Lemma~\ref{lemma liminf
positive} $ u(-\infty)=0$ or $u(+\infty) =0$. Then we may apply
Theorem~\ref{thm exp decay} and deduce the exact asymptotic
behavior of $u$ at either $-\infty$ or $+\infty$ and that $c^*\le
0$ or $c_* \le 0$. Let $u_0$ denote a non-decreasing travelling
wave with speed $c=0$ if $c^*\le 0$ or a non-increasing one if
$c_* \le 0$. Then, by slightly modifying the proof of Theorem~2.1
in \cite{CC} we deduce that for a suitable translation we have
$u^\tau \equiv u_0$. In particular the profile of the travelling
wave $u_0$ is unique.

\fdem

\bigskip

Next we address the issues of non-uniqueness and discontinuities
of solutions when $c=0$. We consider $f$ such that
\begin{align} \label{hyp f not unique 1}
\hbox{$f$ is smooth, $0<f'(0)<1$, $f'(1)<0$ and $f$ is KPP.}
\end{align}
We are interested in the case where  $u-f(u)$ is not monotone, and
for simplicity we shall assume that setting
$$
g(u) = u - f(u)
$$
there exists $0< \alpha<\beta<1$ such that
\begin{align} \label{hyp f not unique 2}
\begin{aligned}
& g'(u) >0 \quad \forall u \in [0,\alpha) \cup (\beta,1]
\\
& g'(u) <0 \quad \forall u \in (\alpha,\beta).
\end{aligned}
\end{align}

\begin{proposition}\label{not unique}
Assume $f$ satisfies \eqref{hyp f not unique 1}, \eqref{hyp f not
unique 2}. Then there exists $J$ such that no solution of
\eqref{mono.eq.tw}-\eqref{mono.eq.twbc+} is continuous, and this
problem admits infinitely many solutions.
\end{proposition}
\dem{Proof.} Let us choose $J \in C^1$, with compact support and
satisfying \eqref{j1} and \eqref{MP}, and such that $c^1 \le 0$.
Then by Corollary~\ref{corollary velocity kpp} we have $c^* = c^1
\le 0$. Thus there exists a monotone travelling wave solution
$u_1$ of \eqref{mono.eq.tw}-\eqref{mono.eq.twbc+} with speed
$c=0$. If \eqref{mono.eq.tw}-\eqref{mono.eq.twbc+} has a
continuous solution $u_2$ , then by Theorem \ref{mono.th.uniq} and
Remark~\ref{one continuous} we have $u_1 \equiv u_2$. Hence $u_1$
is monotone and continuous. Then $J \star u_1$ is monotone which
implies that $u_1 - f(u_1)$ is monotone in $\R$. This is
impossible if $u_1$ is continuous and $u-f(u)$ is not monotone.

For the construction of infinitely many solutions we follow
closely the work of \cite{BFRW}. Since $g'(0)>0$ and $g'(1)>0$
there are $a<b$ such that
$$
\hbox{$g$ is increasing in $[0,a]$, $g$ is increasing in $[b,1]$}
$$
$$
g(a) = g(b) \hbox{ and $g$ is not monotone in $[a,b]$}.
$$
Define
$$
\tilde g(u) = \begin{cases} g(u) & \hbox{if $u \in [0,a]$ or $u\in
[b,1]$ }
\\
g(a) & \hbox{if $u \in [a,b]$}
\end{cases}
$$
Let $g_n:[0,1] \to \R$ be smooth such that $g_n \to g$ uniformly
in $[0,1]$, $g_n \equiv g$ in a neighborhood of 0 and 1, $g_n'>0$
and $u- g_n(u)$ is KPP. Then by Corollary~\ref{corollary velocity
kpp} the problem \eqref{mono.eq.tw}- \eqref{mono.eq.twbc+} with
nonlinearity $f_n = u - g_n(u)$ has critical speed $c^* \le 0$
independent of $n$, and hence there exists a monotone solution
$u_n$
$$
J\star u_n = g_n(u_n), \quad u_n(-\infty) = 0, \quad u_n(+\infty)
=1.
$$
Notice that any solution to this problem is continuous and hence
we may choose
$$
u_n(0) = a.
$$
By Helly's theorem there is a subsequence which converges
pointwise to a solution $u$ of the following problem
$$
J \star u = \tilde g(u) \quad \hbox{in $\R$}.
$$
Remark that $u(0)=a$, and $u(-\infty) = 0$, $u(+ \infty) =1$ by
Lemma~\ref{lema 2.4}. Note that $u$ is continuous in $(-\infty,0]$
since $u\le a$ in $(-\infty,0]$ and $g$ is strictly increasing in
$[0,a]$.

We will show that $u$ has a discontinuity at $0$ and $u(0^+) = b$.
As in \cite{BFRW}, choose $\delta_n >0$ such that $u_n(\delta_n) =
b$. Let $\delta =\liminf \delta_n $ and note that $u \ge b$ in
$(\delta,\infty)$. Let us show that $\delta = 0$. If not, then
$\tilde g(u(x)) = g(a)$ for $x \in (0,\delta)$ and this implies $J
\star u = const$ in $(0,\delta)$. Then for $0<\tau<\delta/2$ we
have $J \star ( u - u(\cdot - \tau) ) \ge 0 $ and vanishes in a
nonempty interval. By the maximum principle $ u \equiv u(\cdot -
\tau) $ and this implies that $u$ is constant, which is a
contradiction. Thus $\delta = 0$ and $u$ has a jump discontinuity
at $0$. Hence $u$ is a solution to \eqref{mono.eq.tw}-
\eqref{mono.eq.twbc+}. We conclude that $u(0^+) = b$ because
$J\star u$ is continuous. \fdem

\bigskip

\noindent{\bf Acknowledgements.} This work has been partly
supported by the Ecos-Conicyt project C05E04, Fondecyt 1050754,
Fondecyt 1050725, Nucleus Millenium P04-069-F, Information and
Randomness and by FONDAP grant for Applied Mathematics, Chile. 
This works also has been partly supported by the Max Planck Institut for mathematics in the science of Leipzig.









\end{document}